\documentclass[10pt]{article}
\usepackage{epsfig,amssymb}

\def\Frac#1#2{\mbox{\Large${\textstyle \frac{#1}{#2}}$}}

\def \ra {{\quad\Rightarrow\quad}}
\def \lr {{\quad\Leftrightarrow\quad}}

\def\R{{\mathbb R}}
\def\N{{\mathbb N}}
\def\A{{\cal A}}
\def\PP{{\cal P}}
\def\QQ{{\cal Q}}
\def\OO{{\cal O}}

\def\wh{\widehat}
\def\e{\epsilon}
\def\la{\lambda}
\def \l {\ell}
\newtheorem{lemma}{Lemma}[section]
\newtheorem{proposition}[lemma]{Proposition}
\newtheorem{corollary}[lemma]{Corollary}
\newtheorem{theorem}[lemma]{Theorem}
\newtheorem{remark}[lemma]{Remark}

\newtheorem{definition}[lemma]{Definition}

\newtheorem{conjecture}[lemma]{Conjecture}
\newtheorem{claim}[lemma]{Claim}
\def\be  {\begin{equation}} 
\def\ee  {\end{equation}}
\def\ba  {\begin{eqnarray}} 
\def\ea  {\end{eqnarray}}
\def\baa {\begin{eqnarray*}} 
\def\eaa {\end{eqnarray*}}
\makeatletter
\@addtoreset{equation}{section}
\makeatother
\def\proof{\medskip\noindent{\bf Proof.} }
\def\qed{\hfill $\Box$}
\newcommand {\lb} {\label}
\newcommand {\rf[1]} {(\ref{#1})}

\usepackage[table]{xcolor}

\def \c {\cellcolor[gray]{0.8}}

\begin{document}

\title{Landau--Kolmogorov inequality revisited}

\author{A.\,Shadrin \\
DAMTP, University of Cambridge, UK}

\date{}
\maketitle


\section{Introduction}


The Landau--Kolmogorov problem consists of finding the upper bound $M_k$ 
for the norm of intermediate derivative $\|f^{(k)}\|$, when the bounds  
$\|f\| \le M_0$ and $\|f^{(n)}\| \le M_n$, 
for the norms of the function and of its higher derivative, 
are given.

Here, we consider the case of a finite interval when $f \in W^n_\infty[-1,1]$
and all the norms are the max-norms, $\|\cdot\| = \|\cdot\|_{L_\infty[-1,1]}$.
Precisely, given $n,k \in \N$ and $\sigma \ge 0$, 
we define the functional class
$$
    W^n_\infty(\sigma) 
:= \{f: f \in W^n_\infty[-1,1]\,,
   \; \|f\| \le 1, \; \|f^{(n)}\| \le \sigma \}\,
$$
and consider the problem of finding the values
\baa
   m_k(x,\sigma) &:=& \sup_{f \in W^n_\infty(\sigma)} |f^{(k)}(x)|\,, 
   \quad x\in [-1,1]\,, \\
   M_k(\sigma) &:=& \sup_{f \in W^n_\infty(\sigma)} \|f^{(k)}\| 
       \;=\; \sup_{x\in[-1,1]} m_k(x,\sigma)\,.
\eaa
Our interest to that particular case is motivated by the fact that 
there are good chances to add this case to a short list of
Landau--Kolmogorov inequalities where a {\it complete} 
solution exists, i.e., a solution that covers {\it all} values
of $n,k\in\N$ (and, for a finite interval, {\it all} values of $\sigma>0$). 
The main guideline in finding out how good these chances are is the 
following conjecture.

\begin{conjecture}[Karlin \cite{ka}]
For all  $n,k\in\N$ and all $\sigma>0$, 
\be \lb{K}
      m_k(1,\sigma) = \sup_{x\in[-1,1]} m_k(x,\sigma)\,.
\ee
\end{conjecture}

If \rf[K] is true for particular set $\{n,k,\sigma\}$, 
then the function $f\in W^{n}_\infty(\sigma)$ that provides extremum 
$M_k(\sigma)$ to the value $\|f^{(k)}\|$ over $W_\infty^n(\sigma)$
is the same as the solution to the pointwise problem 
at the end-point of the interval. The latter solution is  
however known to be a certain Chebyshev or Zolotarev spline 
$Z_n(\cdot,\sigma)$ (which is just a polynomial for small $\sigma$), 
and thus we have a characterization of the extremal function. 

\begin{corollary}
If equality \rf[K] is valid for particular $\{n,k,\sigma\}$, 
then for that set of parameters we have
\be \lb{K_2}
    M_k(\sigma)  
 = \|Z_n^{(k)}(\cdot,\sigma)\| = Z_n^{(k)}(1,\sigma)\,.
\ee
\end{corollary}


So far, Karlin's conjecture has been proved for small $n$ with all $\sigma$, 
and for all $n$ with particular $\sigma$, namely in the following cases:

\begin{center}
\begin{tabular}{lcl}
$n=2$, & all $\sigma$ &  Chui--Smith \cite{cs} ($\sigma \le \sigma_n$),
           Landau \cite{l} ($\sigma > \sigma_n$);  \\
$n=3$, & all $\sigma$, &   Sato \cite{s}, Zvyagintsev--Lepin \cite{zl}; \\
$n=4$, & all $\sigma$, &   Zvyagintsev \cite{z} ($\sigma \le \sigma_n$),
           Naidenov \cite{n} ($\sigma > \sigma_n$); \\
$n \in\N$, & $\sigma = \sigma_n$, &  Eriksson \cite{e}\,.
\end{tabular}
\end{center}

\noindent
Here
$$
    \sigma_n := \|T_n^{(n)}\| = 2^{n-1} n!\,,
$$
where $T_n$ is the Chebyshev polynomial of degree $n$
on the interval $[-1,1]$. 

The value $\sigma=\sigma_n$ serves as a borderline between two types 
of the extremal Zolotarev functions $Z_n(\cdot,\sigma)$: 
if $\sigma \le \sigma_n$, then $Z_n$ is a polynomial of degree $n$, 
while for $\sigma > \sigma_n$ it is a perfect spline of degree $n$ 
with $r$ knots.
There are further borderlines $\sigma_{n,r}$ (with $\sigma_{n,1} := \sigma_n$)
which indicate that the spline $Z_n(\cdot,\sigma)$ has exactly $r$ knots if 
$\sigma_{n,r} < \sigma \le \sigma_{n,r+1}$, but that distinction is hardly 
of any use, since for $n > 3$ there are no reasonable estimates 
for perfect splines even with one knot. 
In this respect, we may apply more or less developed polynomial tools
to tackle the problem for $\sigma \le \sigma_n$, 
and then may try to use polynomial estimates in the spline case, 
when $\sigma > \sigma_n$. 

\medskip
In this paper, we prove Karlin's conjecture in several further subcases. 

\smallskip
1) The first result closes the ``polynomial'' case and proves that, 
for $\sigma \le \sigma_n$, 
the extremum value of the $k$-th derivative of $f \in W^n_\infty(\sigma)$ 
is provided by the corresponding Zolotarev polynomial. 

\begin{theorem} \lb{1}
If
\be \lb{e1}
      n \in \N, \qquad 1 \le k \le n-1\,, \qquad 0 \le \sigma \le \sigma_n\,,
\ee
then Karlin's conjecture \rf[K]-\rf[K_2] is true.
\end{theorem}

2) For the ``spline'' case, we managed to advance only up to the second 
derivative. 

\begin{theorem} \lb{2}
If
$$
      n \in \N, \qquad k=1,2,\qquad \sigma_n < \sigma < \infty,
$$
then Karlin's conjecture \rf[K]-\rf[K_2] is true.
\end{theorem}

\noindent
The further advance depends mostly on improving the lower 
bound for the exact constant $C_{n,k}$ in Landau-Kolmogorov inequality
on the half-line:
$$
    \|f^{(k)}\|_{\R_+} 
\le C_{n,k} \|f\|_{\R_+}^{1-k/n} \|f^{(n)}\|_{\R_+}^{k/n}
$$
The existing lower bounds for $C_{n,k}$,
which are due to Stechkin, are not very satisfactory for general $n$ and 
$k > 2$.

\smallskip
3) However, for small $n$, these bounds 
can be improved, thus leading to one more extension. 

\begin{theorem} \lb{3}
If
\be \lb{n<10}
\begin{array}{ll}
      n = 5,6   & \qquad 1 \le k \le n-2, \\
      n = 7,8,9 & \qquad 1 \le k \le n-3, \\
      n = 10,11 & \qquad 1 \le k \le 6,
\end{array}
\qquad \sigma_n < \sigma < \infty,
\ee
then Karlin's conjecture \rf[K]-\rf[K_2] is true.
\end{theorem}

In all the cases, the proof is based on comparing the upper bound
for the local extrema of the function $m_k(\cdot,\sigma)$ 
with the lower bound for the value $m_k(1,\sigma)$. 
The technique we use is not working for the value $k=n-1$, what explains 
restriction in \rf[n<10]. In \rf[e1], i.e., for $0 \le \sigma <\sigma_n$,
we managed to cover the case $k=n-1$ by different means.

The upper bounds are given in terms of Zolotarev polynomials
and these estimates may be viewed as a generalization to higher derivatives 
of Markov-type results of Schur \cite{schur} and Erd\H{o}s-Szeg\H{o} \cite{es}.
These bounds demonstrate once again, 
if we borrow the words of Shoenberg said about cubic splines, 
``the brave behaviour of Zolotarev polynomials under difficult 
circumstances''.


\section{Main ingredients of the proof}


Karlin's conjecture states that the function $m_k(\cdot,\sigma)$ 
(which is a positive even function)
reaches its maximal value at the end-points of the interval $[-1,1]$. 
To establish this fact it is sufficient to check that, at any point
$x_0$ inside the interval $(-1,1)$ where $m_k(\cdot,\sigma)$ 
takes its local maximum, we have
$$
        m_k(x_0,\sigma) < m_k(1,\sigma)\,.
$$
If $f$ is the function from $W_\infty^n(\sigma)$ that attains a 
locally maximal value $m_k(x_0,\sigma)$, then clearly
$$
    m_k(x_0,\sigma) = |f^{(k)}(x_0)|\,, \qquad f^{(k+1)}(x_0) = 0\,,
$$
and it makes sense to introduce the following quantity:
$$
   m_k^*(x_0,\sigma) 
:= \sup \{ |f^{(k)}(x_0)|:\; f \in W^n_\infty(\sigma),\;
   f^{(k+1)}(x_0) = 0 \},   \qquad x_0 \in [-1,1]\,.
$$
The next statement follows immediately.

\begin{claim}
If, for a given $n,k \in \N$ and $\sigma > 0$, we have
\be \lb{f<m}
     \sup_{x_0\in[-1,1]} m^*_k(x_0,\sigma) \le m_k(1,\sigma)\,,
\ee
then Karlin's conjecture is true.
\end{claim}

In order to verify inequality \rf[f<m], we split it into two 
parts
\be \lb{c1c2}
     m_k^*(x_0,\sigma) \le A(n,k,\sigma)\,,\qquad  
     A(n,k,\sigma) \le m_k(1,\sigma)\,, 
\ee
and then check whether $A \le B$. So, we need two different estimates:

a) a good lower bound for the end-point value  
$m_k(1,\sigma) = \sup \{|f^{(k)}(1)|: f \in W^n_\infty(\sigma)\}$, 

b) a good upper bound for $|f^{(k)}(x_0)|$, where $f$ is from 
$W^n_\infty(\sigma)$ and satisfies $f^{(k+1)}(x_0)=0$.

\noindent
Actually, if $x=x_0$ stays sufficiently far away from the end-points $x=\pm 1$,
then a reasonable upper bound for $|f^{(k)}(x_0)|$ can be established
irrespectively of whether $f^{(k+1)}(x_0)$ vanishes or not. Therefore, 
for the upper bounds for $|f^{(k)}(x)|$, we will consider two cases
$$
   {\rm b}_1) \quad m_k^*(x,\sigma) \le A_{n,k}^*(\sigma), \quad 
      \omega_k < |x_0| \le 1, \qquad
   {\rm b}_2) \quad  m_k(x,\sigma) \le A_{n,k}(\sigma), \quad 
      |x| \le \omega_k < 1,
$$
with an appropriately chosen value $\omega_k$. 

\bigskip
We will distinguish between the cases $\sigma \le \sigma_n$ 
and $\sigma > \sigma_n$.

\bigskip
1)  {\it The case $\sigma \le \sigma_n$.} 

\medskip
1a) {\it Lower estimates for $m_k(1,\sigma)$.}
Clearly, $m_k(1,\sigma)$ is monotoniously increasing with $\sigma$,
therefore, we have the trivial estimate
$$
     m_k(1,\sigma) \ge m_k(1,\sigma_0) = T_{n-1}^{(k)}(1)\,.
$$
However, this estimate is too rough when $k = \OO(n)$, 
so we will use a finer one.

\begin{proposition}
We have
\be \lb{B}
     m_k(1,\sigma) 
\ge  B_{n,k}(\sigma)
:= \left(1 - \Frac{\sigma}{\sigma_n}\right) T_{n-1}^{(k)}(1) 
    +  \Frac{\sigma}{\sigma_n} T_n^{(k)}(1), \qquad
    0 \le \sigma \le \sigma_n.
\ee
\end{proposition}

\proof
Let us show that $m_k(x,\sigma)$ as a function of $\sigma$ is concave. 
For any $x\in[-1,1]$, and for any $\sigma' < \sigma''$, 
let $f_1$ and $f_2$ be the functions such that
$$
   m_k(x,\sigma^{(i)}) = f_i^{(k)}(x), \qquad 
   f_i \in W^n_\infty(\sigma^{(i)}), \qquad i =1,2.
$$
It is clear that, for any $\sigma \in [\sigma',\sigma'']$, 
with $t$ such that $\sigma = (1-t)\sigma' + t \sigma''$, 
the function $f := (1-t) f_1 + t f_2$ belongs to 
$W^n_\infty(\sigma)$, hence we have
$$
     m_k(x,\sigma) \ge f^{(k)}(x) 
= (1 - t) f_1^{(k)}(x) + t f_2^{(k)}(x)
= (1 - t) m_k(x,\sigma') + t  m_k(x,\sigma'')\,.
$$
In particular, with $\sigma_0 := T_{n-1}^{(n)} = 0$ and 
$\sigma_n = T_n^{(n)}$, we have
$$
\textstyle
     m_k(1,\sigma) 
\ge \left(1 - \Frac{\sigma}{\sigma_n}\right) m_k(1,\sigma_0)
    + \Frac{\sigma}{\sigma_n} m_k(1,\sigma_n)\,,
$$
But $m_k(1,\sigma_0) = T_{n-1}^{(k)}(1)$ and 
$m_k(1,\sigma_n) = T_{n}^{(k)}(1)$, hence the result.
\qed

\medskip
1b$_1$) {\it Upper estimate for $m_k^*(x_0,\sigma)$.} 
We will use a comparison lemma of the kind similar to the one that was 
used by Matorin \cite{ma} in (actually) proving that 
$m_k(1,\sigma_n) \le T_n^{(k)}(1)$.

\begin{lemma} \lb{comp}
Let $p \in \PP_n[-1,1]$ be a polynomial that satisfies the following 
conditions:
\be \lb{1-3}
  1)\quad p^{(k+1)}(x_0) = 0\,, \qquad
  2)\quad \mbox{$p$ has an $n$-alternance on $[-1,1]$}, \qquad
  3)\quad \|p^{(n)}\| \ge \sigma\,.
\ee
Then, for any $f \in W^n_\infty[-1,1]$ and for any $x_0 \in [-1,1]$ 
such that
$$
     1') \quad f^{(k+1)}(x_0) = 0\,,\qquad
     2') \quad\|f\| \le 1, \qquad 
     3') \quad \|f^{(n)}\| \le \sigma,
$$
we have 
$$
     |f^{(k)}(x_0)| \le |p^{(k)}(x_0)|\,.
$$
\end{lemma}

\proof
Assume the contrary, i.e., that  $f^{(k)}(x_0) = p^{(k)}(x_0)/\gamma$
with some $\gamma$ such that $|\gamma| < 1$. 
Then the function $g := \gamma f$ satisfies
$$
     2'')\quad \|g\| < 1, \qquad 3'')\quad \|g^{(n)}\| < \sigma\,,
$$
and moreover
$$
    1'')\quad g^{(k)}(x_0) = p^{(k)}(x_0), \qquad
    g^{(k+1)}(x_0) = p^{(k+1)}(x_0) = 0.
$$
Consider the difference $h = p - g$. By the $n$-alternation property $(2)$
of $p$, since $\|g\| < 1$, the function $h$ has at least $n-1$ distinct 
zeros on $[-1,1]$, hence $H := h^{(k-1)}$ has at least $n-k$ distinct zeros
strictly inside $(-1,1)$, and by $(1'')$, we also have 
$H'(x_0) = H''(x_0) = 0$. It follows that $H' = h^{(k)}$ has at least 
$n-k+1$ zeros on $[-1,1]$ counting multiplicities, 
therefore 
$$
    \mbox{$h^{(n)}$ has at least one sign change on $[-1,1]$.}
$$
On the other hand, by $(3)$ and $(3'')$ we have $|g^{(n)}(x)| < \sigma$ 
and $|p^{(n)}(x)| \equiv \mbox{const} \ge \sigma$, hence
$|h^{(n)}(x)| = |p^{(n)}(x) - g^{(n)}(x)| > 0$ for all $x \in [-1,1]$, 
a contradiction. 
\qed

\begin{corollary} \lb{comp1}
We have
\be \lb{p}
    m_k^*(x_0,\sigma) \le |p^{(k)}(x_0)|
\ee
where $p$ is any polynomial of degree $n$ 
that satisfies conditions (1)-(3) in \rf[1-3].
\end{corollary}

\bigskip
Let $\{Z_n(\cdot,\theta)\}$ be the family of 
the Zolotarev polynomials parametrized with respect to the value 
of its highest derivative $\theta := Z_n^{(n)}(\cdot,\theta)$ 
(see Sect.\,\ref{zol} for details).
Given $x_0$, our choice for $p$ in \rf[p] is the dilated Zolotarev polynomial 
$Z_n(\cdot,\theta_{x_0})$ such that $Z_n^{(k+1)}(x_0,\theta_{x_0}) = 0$.
An advantage of choosing such a $p$ is that, for $x_0 \in [\omega_k,1]$, 
the value of $p^{(k)}(x_0)$ can be further bounded in terms of the 
single Zolotarev polynomial $Z_n(\cdot,\theta_k)$ such that
$$
    Z_n^{(k+1)}(1,\theta_k) = 0. 
$$
Namely, as we show in Sects.\,\ref{zol}-\ref{shur}, 
$$
   \sup_{x \in [\omega_k,1]} m_k^*(x_0,\sigma)
\le \max \{1,\Frac{\sigma}{\theta_k}\}^{k/n}
    \max \{ T_n^{(k)}(\omega_k),  Z_n^{(k)}(1,\theta_k)\}
$$
In Sects.\,\ref{k}-\ref{n}, we provide the estimates for the values appeared
here on the right-hand side and, thus, arrive at the following statement.

\begin{proposition} 
We have 
\be \lb{A*}
     \sup_{x \in [\omega_k,1]} m_k^*(x_0,\sigma)
 \le A_{n,k}^*(\sigma) := \left\{ \begin{array}{ll}
     T_{n-1}^{(k)}(1), & 0 \le \Frac{\sigma}{\sigma_n} \le \eta_k; \\
     \la_k T_n^{(k)}(1) \left(
         \Frac{1}{\eta_k}\Frac{\sigma}{\sigma_n}\right)^{k/n},
         & \eta_k \le \Frac{\sigma}{\sigma_n} \le 1.
      \end{array} \right.
\ee
where
$$
   \la_k = \frac{1}{k+1}\frac{n-1}{n-1+k}, \qquad
   \eta_k = \frac{n-(k+1)}{2(2n-(k+1)}\,.
$$
\end{proposition}

\medskip
1b$_2$) {\it Upper estimate for $m_k(x,\sigma)$.} 
We use a technique based on the Lagrange interpolation. 
Let $\ell_\Delta \in \PP_{n-1}$ be the polynomial of degree $n-1$ that 
interpolates $f \in W^n_\infty(\sigma)$ on a mesh $\Delta = (t_i)_{i=1}^n$.
From the identity
$
   f^{(k)}(x) = \ell_\Delta^{(k)}(x) + ( f^{(k)}(x) - \ell_\Delta^{(k)}(x))
$
it follows that
$$
    |f^{(k)}(x)| \le \Lambda_k(x)\|f\| + \Omega_k(x) \|f^{(n)}\|\,,
$$
where
$$
   \Lambda_k(x) = \sup_{\|p\|_\Delta=1} |p^{(k)}(x)|, \qquad
   \Omega_k(x) = \sup_{\|f^{(n)}\|=1} |f^{(k)}(x) - \ell_\Delta^{(k)}(x)|\,,
$$
whence
$$
    \sup_{x \in [0,\omega_k]} m_k(x,\sigma) 
\le \sup_{x \in [0,\omega_k]} \Lambda_k(x) 
    + \sup_{x \in [0,\omega_k]} \Omega_k(x) \sigma\,.
$$
In Sect.\,\ref{x}, we prove that calculation of the suprema on the right-hand
side is reduced to computing the largest local maxima of two specific 
polynomials and that leads to the following estimate.

\begin{proposition} 
We have
\be \lb{A}
    \sup_{x \in [0,\omega_k]} m_k(x,\sigma) 
\le A_{n,k}(\sigma) := \Frac{3}{2k+1} T_{n-1}^{(k)}(1) 
     + \Frac{2}{2k+1}\Frac{2(k+1)}{n+k} T_{n}^{(k)}(1) 
       \Frac{\sigma}{\sigma_n}\,.
\ee
\end{proposition}

The latter estimate is not particularly good for $k=1$ and $k=2$,
so for such $k$ we also use another one
\be \lb{AA}
     \sup_{x \in [0,\omega_k]} m_k(x,\sigma_n)
\le \left(\Frac{1}{1-\sin\frac{k+1}{2n}}\right)^k \Frac{1}{2k+1}
    T_n^{(k)}(1)\,.
\ee

1c) {\it Final step.}
The constants in estimates
\rf[B], \rf[A*] and \rf[A] are easy to compare (they are simple functions 
of $t = \sigma/\sigma_n$) and, in Sect.\,\ref{n-2}, we prove that if 
$n \in \N$, $1 \le k \le n-2$ and $0 \le \sigma \le \sigma_n$,
then 
$$
    \max \left(A_{n,k}(\sigma), A_{n,k}^*(\sigma)\right) 
\le B_{n,k}(\sigma)\,,
$$
and that implies 
$$
    m_k(1,\sigma) =  \sup_{x\in [-1,1]} m_k(x,\sigma), \qquad 
    0 \le \sigma \le \sigma_n, \qquad 1 \le k \le  n-2.
$$

\bigskip
2) {\it The case $\sigma > \sigma_n$.}

\medskip
For that case, it is more convenient to reformulate 
the original problem. Namely, instead of considering functions from the class
$$
    W^n_\infty(\sigma) 
:= \{f: f \in W^n_\infty[-1,1], \;
   \|f\|_{[-1,1]} \le 1, \; \|f^{(n)}\|_{[-1,1]} \le \sigma \}\,, \qquad
   \sigma_n < \sigma < \infty,
$$
i.e., functions on a fixed interval $I_1 = [-1,1]$ with increasing norms 
$\|f^{(n)}\|_{[-1,1]} \le \sigma$, we will consider functions from the class
\be \lb{I_s}
    W^n_\infty(I_s) 
:= \{f:  f \in W^n_\infty[-s,1], \;
   \|f\|_{[-s,1]}\le 1, \; \|f^{(n)}\|_{[-s,1]} \le \sigma_n \}\,, \qquad
   2 < |I_s| < \infty,
\ee
i.e., functions with a fixed norm $\|f^{(n)}\|_{[-s,1]} = \sigma_n$ 
on the intervals $I_s := [-s,1]$ of increasing length $|I_s| > |I_1| = 2$. 
The pointwise Landau-Kolmogorov problem consists then of finding the value
$$
     m_k(x,I_s) := \sup_{f\in W^n_\infty(I_s)} |f^{(k)}(x)|\,,
$$
and Karlin's conjecture states that $m_k(x,I_s)$ is maximal at $x=1$.

\bigskip
2a) {\it Lower estimate for $m_k(1,I_s)$.}
Denote by $B_{n,k}^+$ the best constant in the 
Landau-Kolmogorov inequality on the half-line for the normalized functions:
\ba
     B_{n,k}^+ 
&:=& \sup \{|f^{(k)}(1)|: \|f\|_{[-\infty,1]} \le 1, \;
     \|f^{(n)}\|_{[-\infty,1]} \le \sigma_n \} \,. \lb{B+} \\
&=& \sup \{\|f^{(k)}||_{[-\infty,1]}: \|f\|_{[-\infty,1]} \le 1, \;
     \|f^{(n)}\|_{[-\infty,1]} \le \sigma_n \} \nonumber
\ea

\begin{proposition} 
For all $|I_s| > |I_1| = 2$ we have 
\be \lb{m>B}
    m_k(1,I_s) \ge B_{n,k}^+\,.
\ee
\end{proposition}

\proof
Clearly, with $n$ and $\sigma_n$ fixed, the spaces defined in \rf[I_s]
are embedded into each other, namely
$W^n_\infty(I_s) \supset  W^n_\infty(I_t)$ for  $s < t$,
therefore for the suprema $m_k(1,I_s) := \sup |f^{(k)}(1)|$ 
over those spaces, we have the inequalities
$$
   m_k(1,I_s) \ge  m_k(1,I_t), \qquad s < t. 
$$
Letting $t = -\infty$, we obtain \rf[m>B]. 
\qed

\medskip
2b). {\it Upper estimates for $m_k(x,I_s)$ and $m_k^*(x_0,I_s)$.}
Similar arguments show that the upper bounds for $m_k(x,I_s)$ and 
$m_k^*(x,I_s)$ are majorized by those of $m_k(x,I_1)$ and $m_k^*(x,I_1)$,
respectively. Namely, moving the interval $I = [a,b]$ 
of length $|I| = 2$ inside any $I_s$, we see that 
$W^n_\infty(I_s) \subset W^n_\infty(I)$, hence
\baa
       \sup_{x \in [\omega_k,1]} m_k^*(x,I_s)
&\le&  \sup_{x \in [\omega_k,1]} m_k^*(x,I_1)\,\\
       \sup_{x \in [s_0,\omega_k]} m_k(x,I_s)
&\le&  \sup_{x \in [0,\omega_k]} m_k(x,I_1)\,.
\eaa
where $s_0$ is the middle of the interval $[-s,1]$. 
The right-hand sides are equivalent to the values $m_k^{(*)}(x,\sigma_n)$ 
and for those we have the upper estimates \rf[A*]-\rf[A].
 
\begin{proposition} \lb{mm*}
For all $|I_s| > |I_1| = 2$ we have 
\ba
       \sup_{x \in [\omega_k,1]} m_k^*(x,I_s)
 &\le& A_{n,k}^*(\sigma_n)\,, \lb{A*1} \\
       \sup_{x \in [s_0,\omega_k]} m_k(x,I_s)
 &\le&  A_{n,k}(\sigma_n)\,. \lb{A1}
\ea
\end{proposition}

2c) {\it Final step.}
In Sect.\,\ref{T2} we prove that the constants in 
\rf[m>B]-\rf[A1] satisfy the inequality 
$$
    \max \left(A_{n,k}(\sigma_n), A_{n,k}^*(\sigma_n)\right) 
\le B_{n,k}^+\,, \qquad k=1,2\,
$$
and that proves that 
$$
   m_k(1, I_s) =  \sup_{x\in [-s,1]} m_k(x,I_s), \qquad |I_s| \ge  2,
$$
or, equivalently, 
$$
     m_k(1,\sigma) =  \sup_{x\in [-1,1]} m_k(x,\sigma), \qquad 
     \sigma_n < \sigma < \infty.
$$


\section{Zolotarev polynomials} \lb{zol}


Here, we remind some facts about Zolotarev polynomials taking
some extracts from our survey \cite[p.240-242]{s04}. Note that
we use a slightly different parametrization for $Z_n$.

\begin{definition} \rm
A polynomial $Z_n\in\PP_n$ is called Zolotarev polynomial if 
it has at least $n$ equioscillations on $[-1,1]$, i.e. if there
exist $n$ points
$$
    -1 \le \tau_1 < \tau_2 < \cdots < \tau_{n-1} < \tau_n \le 1
$$
such that  
$$
    (-1)^{n-i} Z_n(\tau_i)= \|Z_n\| = 1.  
$$
\end{definition}

There are many Zolotarev polynomials, for example the Chebyshev polynomials 
$T_n$ and $T_{n-1}$ of degree $n$ and $n-1$, 
with $n+1$ and $n$ equioscillation points, respectively. 
One needs one parameter more to get uniqueness. 
We will use parametrization through the value of the $n$-th derivative 
of $Z_n$:
$$
    \|Z_n^{(n)}\| = \theta
\lr Z_n(x) := Z_n(x,\theta) 
 := \frac{\theta}{n!}x^n + \sum_{i=0}^{n-1} a_i(\theta) x^i\,. 
$$
By Chebyshev's result, $\|p^{(n)}\| \le \|T_n^{(n)}\|\,\|p\|$, 
so the range of the parameter is
$$
    -\sigma_n \le \theta \le \sigma_n, \qquad 
     \sigma_n = \|T_n^{(n)}\| = 2^{n-1}\,n!\,.
$$
As $\theta$ traverses the interval $[-\sigma_n,\sigma_n]$, 
Zolotarev polynomials go through the following transformations: 
$$
-T_n(x) \to  - T_n(ax+b) \to Z_n(x,\theta)
\to  T_{n-1}(x) \to Z_n(x,\theta) \to T_n(cx+d) \to T_n(x)\,.
$$

Zolotarev polynomials subdivide into 3 groups depending 
on the stucture of the set $\A := (\tau_i)$ of their alternation points.

\begin{tabular}{l@{\ }p{0.9\textwidth}}
1) & $\A$ contains $n+1$ points: then $Z_n$ is the Chebyshev polynomial 
     $T_n$.\\
2) & $\A$ contains $n$ points but only one of the endpoints:
     then $Z_n$ is a stretched Chebyshev polynomial $T_n(ax+b)$, 
     $|a| < 1$. \\
3) & $\A$ contains $n$ points including both endpoints:
     then $Z_n$ is called a proper Zolotarev polynomial 
     and it is either of degree $n$, or the Chebyshev polynomial 
     $T_{n-1}$ of degree $n-1$.
\end{tabular}

For a proper Zolotarev polynomial $Z_n$, 
besides the interior alternation points $(\tau_i)_{i=2}^{n-1}$, 
there is a point $\beta = \beta(\theta)$ outside $[-1,1]$ where 
its first derivative vanishes. 

V.\,Markov proved that zeros of $Z_n'(\cdot,\theta)$
are monotonically increasing functions of $\theta \in [-\sigma_n,\sigma_n]$, 
with $\beta$ going through the infinity as $\theta$ passes the zero.
It follows that, for any $\theta_1,\theta_2$, zeros of $Z_n'(\cdot,\theta_1)$ 
and $Z_n'(\cdot,\theta_2)$ interlace with each other, hence by the Markov 
interlacing property the same is true for their derivatives of any order.
In particular, the following lemma is true.

\begin{lemma}
Let $(\alpha_i)_{i=1}^{M-1}$ be the zeros of $T_{n-1}^{(m)}$ in increasing
order, and, for any given $\theta$, let $(\tau_i)_{i=1}^M$ be the zeros of 
$Z_n^{(m)}(\cdot,\theta)$. Then, $(\alpha_i)$ and $(\tau_i)$ interlace, i.e.,
$$
    \tau_1 < \alpha_1 < \tau_2 < \alpha_2 < \tau_3 < 
    \cdots < \alpha_{M-1} < \tau_{M}\,. 
$$
\end{lemma}

Another consequence of the interlacing property is the following observation.

\begin{lemma} \lb{theta}
Let $\omega_k$ be the rightmost zero of $T_n^{(k+1)}$, and let 
$Z_n(\cdot,\theta_k)$ 
be the Zolotarev polynomials whose $(k+1)$st derivative
vanishes at $x=1$, i.e.,
$$
    T_n^{(k+1)}(\omega_k) = 0, \qquad Z_n^{(k+1)}(1,\theta_k) = 0\,.
$$
Further, for a given $x_0 \in (\omega_k,1)$, let $Z_n(\cdot,\theta_{x_0})$
be the Zolotarev polynomial such that
$$
   Z_n^{(k+1)}(x_0,\theta_{x_0}) = 0, \qquad x_0 \in [\omega_k,1]\,.
$$
Then
$$
    |\theta_k| < |\theta_{x_0}| < \sigma_n\,.
$$
\end{lemma}

\proof
According to our parametrization, we have $-T_n(x) = Z_n(x,-\sigma_n)$, and 
as $\theta$ increases from $-\sigma_n$ to $-0$, the rightmost zero
of $Z_n^{(k+1)}(\cdot,\theta)$ increases from $\omega_k$ to $+\infty$,
passing through the value 1 for some $\theta := \theta_k$. 
Therefore
$$
    \omega_k < x_0 < 1 \lr -\sigma_n < \theta_{x_0} < \theta_k.
$$

2) Here we give some upper estimates for the values $T_n^{(k)}(\omega_k)$
relative to the value $T_n^{(k)}(1)$. 
The estimates for $T_n^{(k)}(\omega_k)$ has been given on several 
occasions, we summarize what we need in the following statement.

\begin{lemma}
Let $\omega_k := \omega_{n,k}$ be the rightmost zero of $T_n^{(k+1)}$. Then
\be \lb{Tk}
\begin{array}{ccccl}
  1) & |T_n^{(k)}(\omega_k)| &\le& \Frac{1}{2k+1}\, T_n^{(k)}(1), 
        & n \in\N, \quad 1 \le k \le n-1; \\[1ex]
  2) & |T_n'(\omega_1)|  &\le& \Frac{1}{4}\, T_n'(1),   & n \ge 5; \\[1ex]
  3) & |T_n''(\omega_2)| &\le& \Frac{8}{55}\, T_n''(1), & n \ge 10.
\end{array}
\ee
\end{lemma}

\proof
The first inequality was proved by Eriksson \cite{e} who actually derived
a stronger estimate:
$$
    |T_n^{(k)}(\omega_k)| 
\le \Frac{F_k(\omega_k)}{2k+1} \, T_n^{(k)}(1),
$$
where
$$
     F_k(x) := \frac{2(1+x)^2}{(2k+5)x + 2} \le 1, \qquad x\in[0,1].
$$
The second inequality is due to Erd\"os--Szeg\"o \cite[p.464]{es}. 
To derive the third one, we note that the function $F_k(\cdot)$ 
has the single minimum at $x_* = \frac{2k+1}{2k+5} = \frac{5}{9}$, 
therefore, if  $x_* < \omega_2 < 1$, then 
\be \lb{F}
     F_2(\omega_2) < F_2(1) = \Frac{8}{11}.
\ee
But $\omega_2$ is the largest zero of the third derivative of $T_n$, 
therefore it is greater than the third largest zero of $T_n'$, i.e.,
$\omega_2 > \cos\frac{3\pi}{n}$, so \rf[F] is valid
if $\cos\frac{3\pi}{n} \ge \frac{5}{9}$, 
and the latter holds for $n \ge 10$.
\qed

\begin{corollary}
We have
\be \lb{Tx} 
    \max_{x \in [0,\omega_{k-1}]} |T_n^{(k)}(x)|
\le \Frac{1}{2k+1} T_n^{(k)}(1)
\ee
\end{corollary}

\proof
The values of local maxima of $|T_n^{(k)}(\xi_i)|$ increase 
with $|\xi_i|$, and since $\omega_k = \max_i |\xi_i|$, we have
$$
    \max_{x \in [0,\omega_{k}]} |T_n^{(k)}(x)|
\le |T_n^{(k)}(\omega_k)|
\le \Frac{1}{2k+1} T_n^{(k)}(1)\,
$$
On the interval $[\omega_k, \omega_{k-1}]$ the value $|T_n^{(k)}(x)|$
decreases monotonically from the rightmost maximum $T_n^{(k)}(\omega_k)$ 
to the rightmost zero $T_n^{(k)}(\omega_{k-1}) = 0$, hence the inequality
for such $x$. 
\qed


\section{A generalization of Erd\"os--Szeg\'o result} \lb{shur}


By $\QQ_n$ we denote the unit ball in the space $\PP_n$, i.e., the set
of polynomials $p\in\PP_n$ such that $\|p\|\le 1$.
According to the well-known Markov inequality
$$
    \sup_{p \in \QQ_n} |p'(x)| \le n^2, \quad x \in[-1,1]\,,
$$
and equality is attained at $x=1$ for $p = T_n$. 

In 1913, Schur \cite{schur} considered the problem of finding 
the maximum of $|p'(x_0)|$ under additional assumption that 
$p''(x_0)=0$. Let $\QQ_n^k(x_0)$ be the unit ball of polynomials 
such that $p^{(k+1)}(x_0) = 0$. Shur proved that 
\be \lb{1/2}
    \sup_{ p \in \QQ_n^1(x_0)} |p'(x_0)| < \frac{1}{2}\, n^2\,.
\ee
Moreover, he showed that if $\la_{n}$ is the least constant in front 
of $n^2$, then, for $\la_{\infty} := \lim \sup_{n\to\infty} \la_{n}$, 
we have
$$
    0.217\cdots \le \la_{\infty} \le 0.465\cdots\,.
$$

In 1942, Erd\H{o}s and Szeg\H{o} \cite{es} refined Shur's result by showing
that the limit $\la_{\infty} = \lim_{n\to\infty} \la_{n}$ exists 
and it is equal to
\be
     \la_{\infty} = \kappa^{-2}(1-E/K)^2 = 0.3124\cdots
\ee
where $E,K$ are the complete elliptic integrals associated with the 
modulus $\kappa$. (They did not improve the uniform bound \rf[1/2] though.)

They also showed that, for any $x_0 \in [-1,1]$, the supremum of $|p'(x_0|$
is attained when $p$ is a Zolotarev polynomial $Z_n(\cdot,\theta)$,
and that the maximum over $x_0$ is attained at $x_0 = 1$ for $n \ge 4$, 
and at $x_0=0$ for $n=3$. 

In this section, we generalize these results to the 
derivatives of order $k \ge 2$.

\medskip
Denote by
$$
    \mu_k(x) := \max_{p \in \QQ_n} |p^{(k)}(x)|, \qquad x \in [-1,1],
$$
the best constant in the pointwise Markov inequality, and by
$$
    \mu_k^*(x_0) := \max_{p \in \QQ_n^k(x_0)} |p^{(k)}(x_0)|
    \qquad x_0 \in [-1,1],
$$
the best constant in the pointwise Schur-type inequality. It is clear
that 
$$
    \mu_k^*(x_0) \le \mu_k(x_0), \qquad x_0 \in [-1,1]\,,
$$
and that equality occurs only if $\mu_k'(x_0) = 0$, i.e. if $x_0$ is a point 
of local extremum (maximum or minimum) of the function $\mu_k(\cdot)$ 
inside $(-1,1)$. 

The next two lemmas are straightfroward extensions of the arguments 
given in \cite[pp.461-462]{es}, from $k=1$ to $k \ge 2$.

\begin{lemma} \lb{x_0}
For any $\theta$, if $Z_n^{(k+1)}(x_0,\theta) = 0$, then 
\be \lb{z}
    \mu_k^*(x_0) = Z_n^{(k)}(x_0,\theta).
\ee
Conversely, for any $x_0 \in [-1,1]$, with some $\theta = \theta_{x_0}$
there is a polynomial $Z_n(\cdot,\theta)$ such that \rf[z] is true.
\end{lemma}

\begin{lemma} \lb{x1}
Let $x_0$ be a point such that
$$
     \mu_k^*(x_0) < \mu_k(x_0)\quad {\rm and} \quad x_0 \ne \pm 1\,.
$$
Then, for small $\delta > 0$, there is a point 
$x_1 \in [x_0-\delta,x_0+\delta]$, such that
$$
     \mu_k^*(x_0) < \mu_k^*(x_1)\,.
$$
\end{lemma}

\proof
Let $\mu_k^*(x_0) = Z_n^{(k)}(x_0)$, where $Z_n^{(k+1)}(x_0) = 0$
and let $p \in \QQ_n$ be the polynomial such that
$$
     p^{(k)}(x_0) > Z_n^{(k)}(x_0) > 0\,.
$$
Then the polynomial $q = (1-\e)Z_n + \e p$ satisfies 
$$
     \|q\| \le 1, \qquad q^{(k)}(x_0) > Z_n^{(k)}(x_0) = \mu_k^*(x_0)\,,
$$
and, for small $\e$, its $k$-th derivative has a local maximum in the 
neighbourhood of $x_0$ (because $Z_n^{(k)}$ has). 
Let $x_1$ be the point of that maximum, i.e., $q^{(k+1)}(x_1) = 0$.
Then $q^{(k)}(x_1) > q^{(k)}(x_0)$, and respectively
$$
   \mu_k^*(x_0) < q^{(k)}(x_0) < q^{(k)}(x_1) \le \mu_k^*(x_1)\,,
$$
the latter inequality by definition of $\mu_k^*(\cdot)$.
\qed

\begin{corollary} 
Let $\eta$ be a point of local maximum of the function $\mu_k^*(\cdot)$.
Then
$$
    \mu_k^*(\eta) = \mu_k(\eta).
$$
\end{corollary}

\begin{theorem} \lb{sk}
Let $Z_n(x,\theta_k)$ be the Zolotarev polynomial such that
$$
    Z_n^{(k+1)}(1,\theta_k) = 0.
$$
Then 
$$
    \max_{x_0 \in [-1,1]} \mu_k^*(x_0)
 =  \max\, \{ |T_n^{(k)}(\omega_k)|, |Z_n^{(k)}(1,\theta_k)|\}\,.
$$
\end{theorem}

\proof
Let $\eta_i$ be the points of local maxima of of $\mu_k^*(\cdot)$ inside
the interval $(-1,1)$. Then 
$$
   \max_{x_0 \in [-1,1]} \mu_k^*(x_0)
= \max\,\{\mu_k^*(\eta_i), \mu_k^*(1)\}
$$
The corollary shows that, inside $(-1,1)$, the local maxima 
of $\mu_k^*(\cdot)$ coincide with the extrema (maxima or minima) 
of $\mu_k(\cdot)$. On the other hand, V.\,Markov proved that the local 
maxima of $\mu_k(\cdot)$ coincide with those of $|T_n^{(k)}|$. Hence
$$
   \max_{x_0 \in [-1,1]} \mu_k^*(x_0)
 = \max\,\{ |T_n^{(k)}(\xi_i)|, \mu_k^*(1)\}\,,\quad 
   \mbox{where}\quad T_n^{(k+1)}(\xi_i) = 0\,.
$$
Further, it is known that the local maxima of $|T_n^{(k)}|$
are increasing as $|\xi_i|$ increases, i.e,
$$
    \max_i\,|T_n^{(k)}(\xi_i)| = |T_n^{(k)}(\omega_k)|\,,
$$
where $\omega_k$ is the rightmost zero of $T_n^{(k+1)}$. Finally,
by Lemma \ref{x_0}, 
$$
    \mu_k^*(1) = |Z_n^{(k)}(1,\theta_k)|,
$$
and that completes the proof. 
\qed

\begin{theorem} 
Let $Z_n(x,\theta_k)$ be the Zolotarev polynomial such that
$$
    Z_n^{(k+1)}(1,\theta_k) = 0.
$$
Then 
$$
    \max_{x_0 \in [\omega_k,1]} m_k^*(x_0,\sigma)
\le \max \{1, \Frac{\sigma}{\theta_k} \}^{k/n}
    \max\, \{ |T_n^{(k)}(\omega_k)|, |Z_n^{(k)}(1,\theta_k)| \}\,.
$$
\end{theorem}

\proof
According to Corollary \ref{comp1},
$$
    m_k^*(x_0,\sigma) \le |p^{(k)}(x_0)|,
$$
where $p$ is any polynomial of degree $n$ such that
$$
  1)\quad p^{(k+1)}(x_0) = 0\,, \qquad
  2)\quad \mbox{$p$ has an $n$-alternance in $[-1,1]$}, \qquad
  3)\quad \|p^{(n)}\| \ge \sigma\,.
$$
We take $p$ as a dilated Zolotarev polynomial 
$Z_n(\cdot,\theta_{x_0})$ such that $Z_n^{(k+1)}(x_0,\theta_{x_0}) = 0$.
The latter satisfies conditions (1)-(2), and its highest derivative 
has the value $\theta_{x_0}$. 
So, if $\theta_{x_0} \ge \sigma$, then condtion (3) is fulfilled 
with $p = Z_n(\cdot,\theta_{x_0})$, but if  $\theta_{x_0} < \sigma$, 
then we have to scale $Z_n$ to ensure (3). So we set
$$
    p(x) :=  Z_n(x_0 + \gamma_0^{1/n}(x-x_0), \theta_{x_0}), \qquad 
    \gamma_0 := \max \{1,\Frac{\sigma}{\theta_{x_0}} \}\,,
$$
whence
$$
    m_k^*(x_0,\sigma) \le p^{(k)}(x_0) 
=  \max \{1,(\Frac{\sigma}{\theta_{x_0}})^{k/n}\} 
   Z_n^{(k)}(x_0,\theta_{x_0}).
$$
Finally, 
$$
    \omega_k \le x_0 \le 1
\ra  \left\{ \begin{array}{l}
    1) \quad |Z_n^{(k)}(x_0,\theta_{x_0})| 
\le \max \{ T_n^{(k)}(\omega_k),  Z_n^{(k)}(1,\theta_k)\}\,, \\
    2) \quad |\theta_k| \le |\theta_{x_0}| \le \sigma_n\,,
    \end{array} \right.
$$
where the first inequality us due to Theorem \ref{sk}, 
and the second one is due to Lemma \ref{theta}. 
\qed


\section{Upper estimates for $Z_n^{(k)}(1,\theta_k)$
         and generalization of Schur inequality} \lb{k}


Recall that by Markov's inequality
$$
   \sup_{\|p^{(k)}\|\le 1} |p^{(k)}(x)| \le |T_n^{(k)}(1), \qquad
   x \in [-1,1],
$$
so we will give some upper estimates for the constant $\la_k$ such that
$$
   Z_n^{(k)}(1,\theta_k) \le \la_k T_n^{(k)}(1)
$$
%
We will get those estimates using the following lemma.

\begin{lemma} 
Let $p \in \PP_n$ be any polynomial that satisfies the following 
conditions:
\be \lb{pp}
  1) \quad p^{(k+1)}(1) = 0\,, \qquad
  1)\quad \mbox{$p$ has an $n$-alternance on $[-1,1]$}.
\ee
If $Z_n^{(k+1)}(1,\theta_k) = 0$, then
\be \lb{Zp}
     |Z_n^{(k)}(1,\theta_k)| \le |p^{(k)}(1)|\,.
\ee
\end{lemma}

\proof
The proof is parallel to the proof of Lemma \ref{comp}, since 
$Z_n$ satisfies $\|Z_n\| \le 1$. Assuming the contrary to \rf[Zp], 
we derive that the $n$-th derivative of $h := p - \gamma Z_n$ should
change its sign which is impossible as $h$ is a polynomial of degree $n$
\qed

\bigskip
2a) We will construct several $p$ that satisfy \rf[pp] 
using alternation properties of $T_n$ and 
$T_{n-1}$. We start with the simplest one.

\begin{lemma}
We have
\be \lb{p1}
      |Z_n^{(k)}(1,\theta_k)| \le \frac{1}{k+1} T_n^{(k)}(1)\,.
\ee
\end{lemma}

\proof
Take
$$
    p(x) = T_n(x) - c q(x), \qquad q(x) := (x-1)T_n'(x)\,,\qquad
    p^{(k+1)}(1) := 0,
$$
so that $p$ has an $n$-alternance on $[-\cos\frac{\pi}{n},1]$ for any $c$, 
and where the last equality defines particular 
$c := \frac{T_n^{(k+1)}(1)}{q^{(k+1)}(1)}$. Then
$$
   p^{(k)}(1) 
=  T_n^{(k)}(1) - c q^{(k)}(1)
=  \left(1 - \frac{T_n^{(k+1)}(1)}{q^{(k+1)}(1)}
     \frac{q^{(k)}(1)}{T_n^{(k)}(1)}\right) T_n^{(k)}(1)\,,
$$
and since $q^{(m)}(1) = mT_n^{(m)}(1)$, it follows that
$$
    p^{(k)}(1) 
=  \left(1-\frac{k}{k+1}\right)T_n^{(k)}(1) 
= \frac{1}{k+1}\, T_n^{(k)}(1)\,.
$$
\qed

2b) The next lemma improves the previous estimate for $k = \OO(n)$.

\begin{lemma} \lb{p2}
We have
\ba
      |Z_n^{(k)}(1,\theta_k)| 
&\le& T_{n-1}^{(k)}(1)\,, \lb{p2a} \\
      |Z_n^{(k)}(1,\theta_k)| 
&\le& \frac{1}{k+1}\frac{n-1}{n-1+k} T_{n}^{(k)}(1)\,. \lb{p2b}
\ea
\end{lemma}

\proof
Take
$$
    p(x) = T_{n-1}(x) - c q(x), \qquad q(x) := (x^2-1)T_{n-1}'(x)\,,\qquad
    p^{(k+1)}(1) := 0\,.
$$
Then
$$
   p^{(k)}(1) 
=  T_{n-1}^{(k)}(1) - c q^{(k)}(1)
=  \left(1 - \frac{T_{n-1}^{(k+1)}(1)}{q^{(k+1)}(1)}
     \frac{q^{(k)}(1)}{T_{n-1}^{(k)}(1)}\right) T_{n-1}^{(k)}(1)
=: \wh\la_{n,k} T_{n-1}^{(k)}(1)\,. 
$$
Since $q'(x) = (x^2-1)T''_{n-1}(x) + 2x T'_{n-1}(x)
= x T'_{n-1}(x) + (n-1)^2 T_{n-1}(x)$, we have
$$
   q^{(m)}(1) = T_{n-1}^{(m)}(1) + ((n-1)^2 + (m-1)) T_{n-1}^{(m-1)}(1)\,,
$$
and using 
$$
     T_n^{(k+1)}(1) = \frac{n^2 - k^2}{2k+1} T_n^{(k)}(1), \qquad
     T_n^{(k-1)}(1) = \frac{2k-1}{n^2 - (k-1)^2} T_n^{(k)}(1),
$$
we obtain, after some simplifications,
\baa
    \wh\la_{n,k} 
&=& 1 - \frac{k}{k+1}\frac{(n-1)^2 - k^2}{2(n-1)^2 + (k+1)}
         \frac{2(n-1)^2 + (k-1)}{(n-1)^2 - (k-1)^2} \\
&=& \frac{1}{k+1} + \frac{k}{k+1}
       \frac{4k(n-1)^2 + (k-1)}{((n-1)^2-(k-1)^2)(2(n-1)^2 + (k+1))} \\
&\le& \frac{1}{k+1} + \frac{k}{k+1}\frac{1}{n-k} \le 1\,
\eaa
and that proves the first inequality \rf[p2a]. Using
$$
   T_{n-1}^{(k)}(1) = \gamma T_n^{(k)}(1), \qquad 
   \gamma = \frac{n-1}{n}\frac{n-k}{n-1+k}\,,
$$
we obtain
$$
    \la_{n,k} = \wh\la_{n,k} \gamma 
\le \frac{1}{k+1}\frac{n}{n-k} \gamma
= \frac{1}{k+1}\frac{n-1}{n-1+k}
$$ 
and that proves \rf[p2b].
\qed

\medskip
2c) In the next lemma, we get further improvements for $k=1$ and $k=2$.

\begin{lemma}
We have 
\be \lb{p3}
    Z_n'(1,\theta_1)  \le \frac{1}{3}\, T_n'(1), \qquad
    Z_n''(1,\theta_2)  
\le \frac{3}{\pi^2}\frac{\pi^2-6}{15-\pi^2} <  0.23\, T_n''(1)\,.
\ee
\end{lemma}

\proof
Set $\xi := \cos \frac{\pi}{n}$, and let
$$
    r(x) = T_n(x) - cq(x), \qquad q(x) := (x+1)T_n'(x), \qquad 
    r^{(k+1)}(\xi) := 0\,.
$$
The polynomial $r$ has an $n$-alternance on $[-1,\xi]$, 
so that, after finding $r^{(k)}(\xi)$
we will transform it to the polynomial
$p(x) := r\left(-1 + (x+1)\frac{1+\xi}{2}\right)$,
which has an $n$-alternance on $[-1,1]$ and satisfies
$$
    p^{(k)}(1) =  \left(\frac{1+\xi}{2}\right)^k r^{(k)}(\xi)\,.
$$

Let us find $r^{(k)}(\xi)$. We have 
$$
   r^{(k)}(\xi) 
=  T_n^{(k)}(\xi) - c q^{(k)}(\xi)
=  T_n^{(k)}(\xi) - \frac{q^{(k)}(\xi)}{q^{(k+1)}(\xi)}T_n^{(k+1)}(\xi)\,,
$$
where
$$
   q^{(m)}(\xi) =  (1+\xi)T_n^{(m+1)}(\xi) + m T_n^{(m)}(\xi)\,,
$$
so that setting $a_k := T_n^{(k)}(\xi)$, we obtain
$$
   r^{(k)}(\xi) 
 =  a_k - 
   \frac{(1+\xi)a_{k+1} +  k a_k}{(1+\xi)a_{k+2} + (k+1) a_{k+1} } a_{k+1}\,.
$$
Further, we have
$$
   a_0 = T_n(\xi) = -1, \qquad a_1 = T_n'(\xi) = 0,
$$
and, for $k \ge 2$, the values $a_k$ can be computed from the recurrence
relation
$$
    (\xi^2-1) a_{k+2} + (2k+1)\xi a_{k+1} = (n^2-k^2) a_k\,.
$$
In particular, we find
$$
   a_2 = \frac{n^2}{1-\xi^2}\,, \qquad 
   a_3 = \frac{3 \xi}{1-\xi^2}\,a_2, \qquad
   a_4 = \frac{5\xi}{1-\xi^2}\,a_3 - \frac{n^2-2^2}{1-\xi^2}\,a_2\,.
$$
For $k=1$, this gives
$$
      r'(\xi) 
  =  - \frac{(1+\xi)a_2}{ (1+\xi) a_3 + 2 a_2} a_2
  =  - \frac{n^2}{ 2 + \xi}
\ra
    |p'(\xi)| =  \frac{1+\xi}{2(2 + \xi)} T_n'(1) < \frac{1}{3}T_n'(1)\,.
$$
For $k=2$, we obtain
$$
   r''(\xi)
 = a_2 - \frac{(1+\xi)a_3 + 2 a_2}{ (1+\xi) a_4 + 3 a_3} a_3
\ra p''(1) = c(n,\xi) T_n''(1)\,,
$$
where
$$
   c(n,\xi)
 = \left(\frac{1+\xi}{2}\right)^2 
   \left(\frac{ 6\xi + 3\xi^2}{(2\xi^2 + 9\xi + 4) - n^2(1-\xi^2)} - 1\right) 
   \frac{1}{1-\xi^2} \frac{3}{n^2-1}\,.
$$
One can show that $c(n,\xi) = c(n,\cos\frac{\pi}{n})$ 
is increasing with $n$ to its limit value given in \rf[p3].
\qed

\bigskip
\begin{remark} \rm
We checked two other possibilities to construct $p$. 

1) The option
$$
    p(x) = T_n(x) - cq(x), \qquad q(x) := (x+1)T_n'(x)\,,\qquad 
    p^{(k+1)}(1) := 0\,,
$$
results in
$$
   |p^{(k)}(1)|
 = \frac{1}{2k+1}\frac{4n^2-1}{(2n^2 + (k+1))} T_n^{(k)}(1)\,,
$$
which is slightly worse than \rf[p1].

2) The option
$$
    p(x) = \frac{x-\gamma}{1-\gamma}\,T_{n-1}(x), \qquad
    p^{(k+1)}(1) := 0\,,
$$
is very poor for small $k$, and for large $k = \OO(n)$ 
it is slightly worse than \rf[p2b].
\end{remark}


\section{Lower bound for $Z_n^{(n)}(\cdot,\theta_k)$} \lb{n}


\begin{lemma}
Let $Z_n(x,\theta_k)$ be a Zolotarev polynomial such that 
$$
   Z_n^{(k+1)}(-1,\theta_k) = 0.
$$
Then 
$$
    \theta_k := \|Z_n^{(n)}\| \ge  \eta_{n,k} \sigma_n, \qquad
    \eta_{n,k} := \frac{n-(k+1)}{2(2n-(k+1))}\,.
$$
\end{lemma}

\proof
Set $m = k+1$ and $M = n-m$, and denote by $(\tau_i)_{i=1}^{M}$ the zeros of 
$Z_n^{(m)}$ in increasing order:
$$
    -1 = \tau_1 < \tau_1 < \cdots < \tau_{M} < 1. 
$$
Then 
$$
    Z_n^{(m)}(x) = A (x+1)(x-\tau_2)\cdots(x-\tau_{M})\,,
$$
where
$$
    A 
 = \frac{Z_n^{(m)}(1)}{ 2(1-\tau_2)\cdots(1-\tau_{M})}
 =: \frac{1}{2} \frac{A_1}{A_2}\,,
$$
and respectively
\be \lb{sigma}
    \|Z_n^{(n)}\| = A\,M! = \frac{M!}{2} \frac{A_1}{A_2}\,.
\ee

Let us find lower bounds for the constants $A_1$ and $1/A_2$.

1) Let $(\alpha_i)_{i=1}^{M-1}$ be the zeros of $T_{n-1}^{(m)}$ in increasing
order. They interlace with zeros of $Z_n^{(m)}$, i.e.
$$
    -1 = \tau_1 < \alpha_1 < \tau_2 < \alpha_2 < \tau_3 < 
    \cdots < \alpha_{M-1} < \tau_{M} < 1, 
$$
therefore 
$$
   \frac{1}{A_2} := \frac{1}{(1-\tau_2)\cdots(1-\tau_{M})}
 >  \frac{1}{(1-\alpha_1)\cdots(1-\alpha_{M-1})}
$$
On the other hand,
$$
   T_{n-1}^{(m)}(x) 
= \frac{\|T_{n-1}^{(n-1)}\|}{(M-1)!}(x-\alpha_1)\cdots(x-\alpha_{M-1})
\ra
   T_{n-1}^{(m)}(1) 
= \frac{\|T_{n-1}^{(n-1)}\|}{(M-1)!}(1-\alpha_1)\cdots(1-\alpha_{M-1})\,,
$$
and respectively
\be \lb{A_2}
   \frac{1}{A_2} > \frac{1}{(1-\alpha_1)\cdots(1-\alpha_{M-1})}
 = \frac{1}{(M-1)!}\frac{\|T_{n-1}^{(n-1)}\|}{T_{n-1}^{(m)}(1)}\,.
\ee

2) The lower bound for $A_1$ is provided by 
\be \lb{A_1}
   A_1 := Z_n^{(m)}(1,\theta_k) 
\ge T_{n-1}^{(m)}(1) \frac{\sigma_n - \theta_k}{\sigma_n}
   + T_{n}^{(m)}(1) \frac{\theta_k}{\sigma_n}
= \frac{T_{n-1}^{(m)}(1)}{\sigma_n} 
    \left( (\sigma_n-\theta_k) 
       + \frac{T_n^{(m)}(1)}{T_{n-1}^{(m)}(1)}\theta_k\right).
\ee

3) Combining estimates \rf[sigma]-\rf[A_1], we obtain
$$
    \theta_k
\ge \frac{n-m}{2} \frac{\|T_{n-1}^{(n-1)}\|}{\sigma_n}
    \left( (\sigma_n-\theta_k) 
        + \frac{T_n^{(m)}(1)}{T_{n-1}^{(m)}(1)}\theta_k\right).
$$
From the relations
$$
   \frac{\|T_{n-1}^{(n-1)}\|}{\sigma_n} 
:=  \frac{\|T_{n-1}^{(n-1)}\|}{\|T_n^{(n)}\|} = \frac{1}{2n}, \qquad
   \frac{T_{n}^{(m)}(1)}{T_{n-1}^{(m)}(1)} 
=  \frac{n}{n-1}\frac{n-1+m}{n-m}
> \frac{n+m}{n-m}\,,
$$
it follows that
$$
      \theta_k
 >  \frac{n-m}{4n} \left( \sigma_n - \theta_k
      + \frac{n+m}{n-m} \theta_k \right) \\
 =  \frac{n-m}{4n} \left( \sigma_n
      + \frac{2m}{n-m}\theta_k\right)\,.
$$
So, $(1 - \frac{m}{2n})\theta_k \ge \frac{n-m}{4n} \sigma_n$,
and finally
$$
   \theta_k > \frac{n-m}{2(2n-m)}\, \sigma_n, \qquad m = k+1\,.
$$

\begin{proposition} 
We have 
\be \lb{AA*}
     \sup_{x \in [\omega_k,1]} m_k^*(x_0,\sigma)
 \le A_{n,k}^*(\sigma) := \left\{ \begin{array}{ll}
     T_{n-1}^{(k)}(1), & 0 \le \Frac{\sigma}{\sigma_n} \le \eta_k; \\
     \la_k T_n^{(k)}(1) \left(
         \Frac{1}{\eta_k}\Frac{\sigma}{\sigma_n}\right)^{k/n},
         & \eta_k \le \Frac{\sigma}{\sigma_n} \le 1.
      \end{array} \right.
\ee
where
$$
   \la_k = \frac{1}{k+1}\frac{n-1}{n-1+k}, \qquad
   \eta_k = \frac{n-(k+1)}{2(2n-(k+1)}\,.
$$
\end{proposition}


\section{Upper estimates for $m_k(x)$ for $x \in [0,\omega_k]$} \lb{x}


\begin{lemma}
We have
$$
    m_k(x) 
\le \frac{3}{2k+1} T_{n-1}^{(k)}(1)
      + \frac{2}{2k+1} \frac{2(k+1)}{n+k} T_n^{(k)}(1) 
         \frac{\sigma}{\sigma_n}\,.
$$
\end{lemma}

\proof
For $f \in W^n_\infty(\sigma)$, let $\l\in \PP_n$ be the Lagrange polynomial
of degree $n$ that interpolates $f$ at the points of local extrema of 
$T_{n-1}$ on the interval $[-1,1]$, i.e.
$$
    \l(x) = f(x), \qquad (x^2-1)T_{n-1}'(x) = 0\,.
$$
Then
$$
   f^{(k)}(x) 
= \l^{(k)}(x) + (f^{(k)}(x) - \l^{(k)}(x))
\le D_k(x)\|f\| + \Omega_k(x) \|f^{(n)}\|\,,
$$
where
$$
   D_k(x) := \sup_{\|p_{n-1}\|_*=1} |p_{n-1}^{(k)}(x)|,\qquad
   \Omega_k(x) := \sup_{\|f^{(n)}\|=1} |f^{(k)}(x) - \l^{(k)}(x)|\,.
$$

1) For the first constant, we have the estimate
$$
    D_k(x) \le \max \{U(x), V(x)\},
$$
where $U(x) := |T_{n-1}^{(k)}(x)|$ and
\baa
     V(x) 
&:=& \Big|\frac{1}{k}\,(x^2-1)\,T_{n-1}^{(k+1)}(x) + xT_{n-1}^{(k)}(x)\Big| \\
&\le& \frac{k-1}{k}|T_{n-1}^{(k)}(x)| 
      + \frac{(n-1)^2-(k-1)^2}{k}|T_{n-1}^{(k-1)}(x)|\,.
\eaa
We have
\baa
   U(x) &\le& \frac{1}{2k+1} T_{n-1}^{(k)}(1)\,, \\
      V(x) 
&\le& \frac{k-1}{k}\frac{1}{2k+1} T_{n-1}^{(k)}(1)
      + \frac{(n-1)^2-(k-1)^2}{k}\frac{1}{2k-1}T_{n-1}^{(k-1)}(1) \\
& = & \left(\frac{k-1}{k}\frac{1}{2k+1} + \frac{1}{k}\right) 
        T_{n-1}^{(k)}(1) \\
& = & \frac{3}{2k+1}  T_{n-1}^{(k)}(1)\,.
\eaa

2) For the second constant, we have
$$
   \Omega_k(x) \le \max |\frac{1}{n!}\omega^{(k)}(x)|\,.
$$
where $\omega(x) = c (x^2-1)T_{n-1}'(x)$,
with its leading coefficient equal to one, i.e., 
$c = \frac{1}{2^{n-2}} \frac{1}{n-1}$.
Set 
$$
   q(x) := (x^2-1)T_{n-1}'(x)\,.
$$
Then 
$$
   \Omega_k(x) 
\le \frac{1}{2^{n-2}}\frac{1}{n!}\frac{1}{n-1}\max |q^{(k)}(x)|
 = \frac{2}{\sigma_n}\frac{1}{n-1}\max |q^{(k)}(x)|\,. 
$$
Since $q'(x) = (n-1)^2 T_{n-1}(x) + x T_{n-1}'(x)$, we have 
\baa
   q^{(k)}(x) 
&=&  ((n-1)^2 + (k-1))T_{n-1}^{(k-1)}(x) + x T_{n-1}^{(k)}(x) \\
&\le& \frac{(n-1)^2 + (k-1)}{2k-1} T_{n-1}^{(k-1)}(1) 
         + \frac{1}{2k+1} T_{n-1}^{(k)}(1) \\
&=& \left(\frac{(n-1)^2 + (k-1)}{(n-1)^2 - (k-1)^2} 
         + \frac{1}{2k+1}\right)  T_{n-1}^{(k)}(1)
 = \frac{c_{n,k}}{2k+1} T_n^{(k)}(1)\,,
\eaa
where
$$
   c_{n,k} 
= \frac{2(k+1)(n-1)^2 + (k+2)(k-1)}{(n-1+k)(n-1+(k-1))}\frac{n-1}{n}
\le 2(k+1)\frac{n-1}{n-1+k}\frac{n-1}{n}
\le 2(k+1)\frac{n-1}{n+k}\,.
$$
Thus
$$
   \Omega_k(x) 
\le \frac{2}{2k+1} \frac{2(k+1)}{n+k} \frac{1}{\sigma_n} T_n^{(k)}(1)\,.
$$
\qed

\begin{corollary}
We have
$$
    m_k(x,\sigma_n) \le \frac{3}{2k+1}T_n^{(k)}(1)\,, \quad k \ge 2\,.
$$
\end{corollary}

\proof
We have
$$
    m_k(x,\sigma_n) \le \alpha_{n,k} T_n^{(k)}(1),
$$
where
\baa
   \alpha_{n,k} 
&=& \frac{3}{2k+1}\frac{n-1}{n}\frac{n-k}{n-1+k} 
    + \frac{2}{2k+1}\frac{2(k+1)}{n+k}
\le \frac{3}{2k+1}\frac{n-k}{n+k} + \frac{2}{2k+1}\frac{2(k+1)}{n+k} \\
&=& \frac{3}{2k+1}\frac{3n + k + 4}{3n+3k} \le \frac{3}{2k+1}\,.
\eaa

\begin{lemma}
We have
$$
    \max_{x \in [0,\omega_k]} m_k(x,\sigma_n) 
\le \frac{1}{(1-\delta_k/2)^k} T_n^{(k)}(\omega_k)\,,
$$
where $\delta_k$ is the maximal distance between two consecutive zeros of
$T_n^{(k+1)}$.
\end{lemma}

\proof 
We will use the following estimate.
Let $f\in W^n_\infty(\sigma_n)$, i.e.,
$\|f\| \le 1$ and $\|f^{(n)}\| \le \|T_n^{(n)}\|$. Then
$$
   T_n^{(k+1)}(\xi_i) = 0 \ra 
   |f^{(k)}(\xi_i)| \le  |T_n^{(k)}(\xi_i)| \le T_n^{(k)}(\omega_k)\,.
$$
Let $(\xi_i)$ be the zeros of $T_n^{(k+1)}$, and let 
$\delta_k = \max_i |\xi_i - \xi_{i+1}|$.
Set 
$$
    \wh T_n(x) =  T_n(\gamma x), \qquad 
    \gamma = \frac{1}{(1-\delta_k/2)} > 1.
$$
Then 
$$
   \wh T_n^{(k+1)}(\xi) = 0 \ra 
   |f^{(k)}(\xi)| 
\le  |\wh T_n^{(k)}(\xi)| \le \gamma^k T_n^{(k)}(\omega_k)\,.
$$

\begin{corollary}
We have
$$
     \max_{x\in[0,\omega_k]} m_k(x,\sigma_n)
\le  \left(\frac{1}{1-\sin\frac{\pi(k+1)}{2n}}\right)^k 
     T_n^{(k)}(\omega_k)\,.
$$
\end{corollary}

\proof
Since $(\cos\frac{\pi i}{n})$ are zeros of $T_n'$, the zeros $\xi_i$ 
of $T_n^{(k+1)}$ are located in the intervals
$\cos\frac{\pi(i+k)}{n} < \xi_i < \cos\frac{\pi i}{n}$, and 
for the distance between two consecutive $\xi_i$ we have
$$
    \delta_k 
 =  \max_i |\xi_i - \xi_{i+1}| 
\le \max_i \left|\cos\frac{\pi i}{n} - \cos\frac{\pi(i+(k+1))}{n}\right|
\le 2 \sin\frac{\pi(k+1)}{2n}\,.
$$

\begin{corollary}
We have
\be \lb{m1}
   \max_{x\in[0,\omega_1]} m_1(x,\sigma_n) \le \frac{1}{2} T_n'(1)\,.
\ee
\end{corollary}

\proof
a) For $n = 4$, we have
$$
    T_4(x) = 8x^4 - 8x^2 + 1, \qquad 
    T_4'(x) = 16(2x^3-x), \qquad T_4''(x) = 16(6x^2-1),
$$
so that
$$
    \omega_1 = 1/\sqrt{6}, \qquad
    \delta_1/2 = 1/\sqrt{6}, \qquad
    T_4'(\omega_1) = 32/3\sqrt{6} = 2/3\sqrt{6} T_n'(1)\,,
$$
hence
$$
   \alpha_4 =  \frac{1}{1-1/\sqrt{6}}\frac{2}{3\sqrt{6}} 
< 0.46 \le 0.5\,.
$$

b) For $n = 5$, we have
$$
    T_t(x) = 16x^5 - 20x^3 +x, \qquad 
    T_5'(x) = 5(16x^4-12x^2+1), \qquad T_5''(x) = 40x(8x^3-3),
$$
so that
$$
   \omega_1 = \sqrt{\frac{3}{8}}, \qquad
   \delta_1/2 = \frac{1}{2}\sqrt{\frac{3}{8}}, \qquad
   T_5'(\omega_1) = \frac{25}{4} = \frac{1}{4} T_5'(1)\,,
$$
and
$$
   \alpha_5 =  \frac{1}{1-\frac{1}{2}\sqrt{\frac{3}{8}}}\frac{1}{4}
< 0.361 \le 1/2\,.
$$

c) For $n \ge 6$, we have $T_n'(\omega_1) \le \frac{1}{4}T_n'(1)$, 
hence
$$
    \alpha_n \le \frac{1}{1-\sin\frac{\pi}{n}}\frac{1}{4} 
\le 1/2\,.
$$
\qed


\section{Proof of Theorem \ref{1}, the case $k \le n-2$ \lb{n-2}}


\begin{theorem}
We have
$$
     \max_{x_0\in[\omega_k,1]} m_k^*(x_0,\sigma) \le  m_k(1,\sigma), 
     \qquad 0 \le \sigma \le \sigma_n.
$$
\end{theorem}

\proof
1) {\it The case $\sigma \le \theta_k$.} By Lemma \ref{p2}, we have
$$
      m_k^*(x_0,\sigma) \le T_{n-1}^{(k)}(1)\,,
$$
while
$$
   m_k(1,\sigma) > m_k(1,\sigma_0) =  T_{n-1}^{(k)}(1).
$$

\bigskip
2) {\it The case $\sigma > \theta_k$.}
In this case 
$$
     m_k^*(x_0,\sigma) 
\le \frac{1}{k+1}\frac{n-1}{n-1+k}
    \left( \frac{\sigma}{\theta_k}\right)^{k/n} T_n^{(k)}(1)
=  \gamma \left( \frac{t}{\alpha}\right)^{k/n} T_n^{(k)}(1)\,,
$$
and
$$
    m_k(1,\sigma) 
\ge  (1 - t) T_{n-1}^{(k)}(1) + t T_{n}^{(k)}(1)
= \left(\beta(1 - t) + t \right) T_n^{(k)}(1)\,,
$$
where
$$ 
    \alpha := \frac{n-(k+1)}{2(2n-(k+1))}\,, \qquad
    \beta 
:= \frac{T_{n-1}^{(k)}(1)}{T_n^{(k)}(1)} 
 = \frac{n-1}{n}\frac{n-k}{n-1+k}\,, \qquad
    t := \sigma/\sigma_n\,.
$$
So, we need to prove that
$$
     f(t) := \gamma \left(\frac{t}{\alpha}\right)^{k/n}
\le  \beta(1-t) + t =: g(t)\,,\qquad
     t \in [\alpha,1]\,.
$$
The function $f$ is concave, therefore it is bounded from above by its 
tangent $\ell$ at $t=2\alpha$, i.e.
$$
    f(t) \le \ell(t) 
= \gamma 2^{k/n}\left( 1 + \frac{k}{n}\frac{t-2\alpha}{2\alpha} \right)\,.
$$
So, we are done, once we prove that 
$$
     \ell(t) \le g(t) \quad \mbox{on} [\alpha,1].
$$ 
Both functions are straight lines, so we need to check this inequality 
only at the end-points. 

1) At $t=\alpha$, we have 
$$
   \ell(\alpha) 
 =  \gamma 2^{k/n}\left(1 - \frac{k}{2n}\right) 
\le \gamma \left(1 + \frac{k}{n}\right)\,, \qquad
   g(\alpha) \ge g(0) = \beta\,.
$$
So, we need the inequality
$$
    \gamma \frac{n+k}{n} \le \beta \lr 
     \frac{1}{k+1}\frac{n-1}{n-1+k} \frac{n+k}{n}
 \le \frac{n-1}{n}\frac{n-k}{n-1+k} \lr
     \frac{n+k}{k+1} \le n-k\,,
$$
amd the latter is valid for $k \le n-2$.

2)  At $t=1$, we have $g(1) = 1$, while
$$
    \ell(1) 
= \gamma 2^{k/n}\left( 1 + \frac{k}{n}\frac{1-2\alpha}{2\alpha} \right)
= \frac{1}{k+1}\frac{n-1}{n-1+k} 2^{k/n}
     \left( 1 + \frac{k}{n}\frac{n}{(n-(k+1))} \right)\,.
$$
Expression in the parenthesis is less than $1+k$, so 
$$
    \ell(1) 
\le 2^{k/n} \frac{n-1}{n-1+k}
\le \frac{n+k}{n}  \frac{n-1}{n-1+k} < 1.
$$

\begin{theorem}
We have
$$
     \max_{x\in[0,\omega_k]} m_k(x,\sigma) \le  m_k(1,\sigma), 
     \qquad 0 \le \sigma \le \sigma_n.
$$
\end{theorem}

\proof
1) For $k \ge 2$, we use the estimates
$$
    m_k(x,\sigma)
\le \frac{3}{2k+1}T_{n-1}^{(k)}(1) 
    + \frac{2}{2k+1}\frac{2(k+1)}{n+k}\,t\,T_n^{(k)}(1)
=: \ell_1(t)\,.
$$
and
$$
    m_k(1,\sigma)
\ge (1-t)T_{n-1}^{(k)}(1) + t T_n^{(k)}(1)
 =: \ell_2(t), \qquad 
$$
To prove that $\ell_1(t) \le \ell_2(t)$ it is sufficient to compare their
values at the end-points:
\baa
   \ell_1(0) = \frac{3}{2k+1}T_{n-1}^{(k)}(1) 
\le T_{n-1}^{(k)}(1) = \ell_2(0)\,, \\
   \ell_1(1) \le \frac{5}{2k+1}T_{n}^{(k)}(1) 
\le T_{n}^{(k)}(1) = \ell_2(1)\,.
\eaa

2) For $k=1$, we use the following estimates:
$$
   m_k(1,\sigma) \ge m_k(1,\sigma_0) 
= T'_{n-1}(1) = \frac{(n-1)^2}{n^2} T'_n(1) \ge \frac{9}{16} T_n'(1)\,,
$$
and
$$
   m_k(x,\sigma) \le m_k(x,\sigma_n) \le \frac{1}{2} T_n'(1)\,.
$$


\section{Proof of Theorem \ref{1}: the case $k=n-1$} \lb{n-1}


Here we cover the case $k=n-1$ for $0 \le \sigma \le \sigma_n$. 

\begin{theorem}
We have
$$
   m_{n-1}(x,\sigma) \le m_{n-1}(1,\sigma) = Z_n^{(n-1)}(1,\sigma), \qquad 
   0 \le \sigma \le \sigma_n\,.
$$
\end{theorem}

\proof
For $f \in W^n_\infty(\sigma)$, let $\l\in \PP_{n-1}$ be the Lagrange 
polynomial of degree $n-1$ that interpolates $f$ at the points of 
local extrema of $Z_n(\cdot,\sigma)$ on the interval $[-1,1]$, i.e.
$$
    \l(\tau_i,\sigma) = f(\tau_i), \qquad 
    -1 = \tau_0 < \tau_1 < \cdots < \tau_{n-2} < \tau_{n-1} = 1.
$$
Then
$$
   f^{(n-1)}(x) 
= \l^{(n-1)}(x,\sigma) + (f^{(n-1)}(x) - \l^{(n-1)}(x,\sigma))
\le D_{n-1}(x,\sigma)\|f\| + \Omega_{n-1}(x,\sigma) \|f^{(n)}\|\,,
$$
where
$$
   D_{n-1}(x,\sigma) := \sup_{\|p_{n-1}\|_*=1} |p_{n-1}^{(n-1)}(x)|,\qquad
   \Omega_{n-1}(x,\sigma) 
:= \sup_{\|f^{(n)}\|=1} |f^{(n-1)}(x) - \l^{(n-1)}(x,\sigma)|\,.
$$
Therefore, 
\be \lb{m}
   m_{n-1}(x,\sigma) \le D_{n-1}(x,\sigma) + \Omega_{n-1}(x,\sigma) \sigma\,.
\ee

1) It is known that the extremum value $D_{n-1}(x,\sigma)$ 
(which is a constant, since $p^{(n-1)} \equiv {\rm const}$)
is attained by the polynomial $p \in \PP_{n-1}$ such that 
\be \lb{p*}
   p(\tau_i,\sigma) = (-1)^i, \qquad i = 0,\ldots,n-1\,.
\ee
It is easy to see that, with 
$$
   \omega(x,\sigma) := \prod (x-\tau_i),
$$
we have
$$
   p(x) = Z_n(x,\sigma) - \frac{\sigma}{n!} \omega(x,\sigma)\,.
$$
Indeed, \rf[p*] is clearly fulfilled, and $p$ is of degree $n-1$
because the leading coefficients of both polynomials on the right-hand side
are equal to $\sigma/n!$. Therefore
\be \lb{D}
   D_{n-1}(x,\sigma) =  p^{(n-1)}(1,\sigma) 
 = Z_n^{(n-1)}(1,\sigma) - \frac{\sigma}{n!} \omega^{(n-1)}(1,\sigma) > 0.
\ee

2) For $\Omega_{n-1}(x,\sigma)$ we show below that 
\be \lb{Om}
    \Omega_{n-1}(x,\sigma) 
\le \Omega_{n-1}(1,\sigma) = \frac{1}{n!} \omega^{(n-1)}(1,\sigma)\,.
\ee

Thus, from \rf[m]-\rf[Om], we obtain
$$
   m_{n-1}(x,\sigma) 
\le |Z_n^{(n-1)}(1,\sigma) - \frac{\sigma}{n!} \omega^{(n-1)}(1,\sigma)|
 + |\frac{\sigma}{n!} \omega^{(n-1)}(1,\sigma)|
=  Z_n^{(n-1)}(1,\sigma)\,,
$$
and theorem is proved.
\qed

\begin{lemma}
We have
\be \lb{Om1}
    \Omega_{n-1}(x,\sigma) 
\le \Omega_{n-1}(1,\sigma) = \frac{1}{n!} \omega^{(n-1)}(1,\sigma)\,.
\ee
\end{lemma}

\proof
For $\Omega_{n-1}(x,\sigma)$ we have the convex majorant
$$
    \Omega_{n-1}(x,\sigma) 
\le \Omega_{n-1}^*(x,\sigma)
 =  \frac{1}{n}\sum_{i=0}^{n-1} |x - \tau_i(\sigma)|\,,
$$
so that 
$$
    \Omega_{n-1}(x,\sigma) 
\le \max\{ \Omega_{n-1}^*(0,\sigma), \Omega_{n-1}^*(1,\sigma)\}
$$
We note that
$$
   \Omega_{n-1}^*(1,\sigma)
 =  1 - \frac{1}{n}\sum \tau_i(\sigma) 
 = \frac{1}{n!} |\omega^{(n-1)}(1,\sigma)|
 = \Omega_{n-1}(1,\sigma)\,,
$$
so we need to prove that 
$$
     c_1(\sigma) := \frac{1}{n}\sum_{i=1}^{n} |\tau_i(\sigma)|
\le  1 - \frac{1}{n}\sum_{i=1}^n \tau_i(\sigma) =: c_2(\sigma)\,.
$$
For large $n$, this inequality is self-evident because the alternation
points $\tau_i(\sigma)$ are spread sufficiently uniform in the interval
$[-1,1]$, therefore $c_1(\sigma) < 1$ while $c_2 \to 1$. 
But we need it for all $n \ge 2$.

We will use the monotonicity property of $\tau_i(\sigma)$ as functions
of $\sigma$. We have
\be \lb{tau1}
     \tau_i(\sigma_0) \le \tau_i(\sigma) \le \tau_i(\sigma_n)\,.
\ee
Here, $\tau_i(\sigma_0)$ are zeros of $(x^2-1)T_{n-1}'(x)$ and 
$\tau_i(\sigma_n)$ are zeros of $(x-1)T_n'(x)$, therefore
$$
   \cos\frac{\pi((n-i)}{n-1} 
\le \tau_i(\sigma) \le \cos\frac{\pi(n-i)}{n}\,, 
   \quad i = 1,\ldots,n-1, \qquad \tau_{n}(\sigma) = 1.
$$
It follows that
$$
     c_2(\sigma)
 = 1 - \frac{1}{n} \sum \tau_i(\sigma) 
 \ge 1 - \frac{1}{n} \sum \tau_i(\sigma_n) 
= 1 - \frac{1}{n}\,.
$$
On the other hand, with $m=\lfloor\frac{n}{2} \rfloor$, 
$$
      \sum |\tau_i(\sigma)|
\le \sum_{i=1}^m |\tau_i(\sigma_0)|
       + \sum_{i=m+1}^{n} |\tau_i(\sigma_n)|
 = \sum_{i=0}^{m-1} \cos\frac{\pi i}{n-1}
    + \sum_{i=0}^{m-1} \cos\frac{\pi i}{n}
\le 1 + \frac{1}{\sin\frac{\pi}{2n}}\,,
$$
where we used the inequality
$$
   \sum_{i=0}^{m-1} \cos ix 
=  \frac{1}{2} + \left(\frac{1}{2} + \sum_{i=1}^{m-1} \cos ix \right)
=  \frac{1}{2} + \frac{\sin(m-\frac{1}{2})x}{2\sin\frac{1}{2}x}
\le  \frac{1}{2} + \frac{1}{2\sin\frac{\pi}{2n}}, \qquad 
    x \in \{\frac{\pi}{n},\frac{\pi}{n-1}\}\,.
$$

a) For $n \ge 6$ we have
$$
    c_1(\sigma)
\le \frac{1}{n} + \frac{1}{n\sin\frac{\pi}{2n}}
\le \frac{1}{6} + \frac{1}{6\sin\frac{\pi}{12}}
 = 0.81 < \frac{5}{6} < 1 - \frac{1}{n} \le c_2(\sigma)\,.
$$

b) For $n=5$,
$$
    c_1(\sigma)
\le \frac{1}{5}\left( 1 + \cos\frac{\pi}{4} 
    + 1 + \cos\frac{\pi}{5} + \cos\frac{2\pi}{5} \right)
 = 0.76 < \frac{4}{5} = \Omega_{n-1}(1,\sigma)\,.
$$

c) For $n=3$ and $n=4$, we cannot obtain the inequality
$c_1(\sigma) \le c_2(\sigma)$ through the estimates \rf[tau1].
In these cases we split the interval $[\sigma_0,\sigma_n]$ into two parts: 
$$
    1) \quad \tau_i(\sigma_0) \le \tau_i(\sigma) \le \tau_i(\wh\sigma_n)\,,
    \quad \sigma \in [\sigma_0,\wh\sigma_n]; \qquad
    2) \quad \sigma \in [\wh\sigma_n,\sigma_n]\,,
$$
where the second interval conatins $\sigma$ such that 
$Z_n(\cdot,\sigma)$ are the Chebyshev polynomials
stretched from the interval $[-\cos\frac{\pi}{n},1]$ to a slightly 
larger interval $[-\cos \phi,1]$ up to $[-1,1]$, i.e.
$$
   Z_n(x,\sigma) = T_n(1 + s(x-1))\,, \qquad
   s \in [s_n,1], \qquad
   s_n := \frac{1+\cos\frac{\pi}{n}}{2} = \cos^2\frac{\pi}{2n}\,.
$$
The alternation points of such $Z_n$ are given by
$$
    \tau_i(\sigma) 
= (1 + t)\cos\frac{(n-i)\pi}{n} - t, \qquad
  t \in[0,t_n], \qquad
  t_n = \tan^2\frac{\pi}{2n}\,.
$$

c$_1$) Consider first the case $\sigma \in [0,\wh\sigma_n]$. 

For $n=2$, 
$$
   \tau_1(\sigma) = -1, \qquad \tau_2(\sigma) = 1\,.
$$

For $n=3$, we have
$$
   \tau_1(\sigma) = -1, \quad
   0 \le \tau_2(\sigma) \le \frac{1}{3}, \qquad
   \tau_3(\sigma) = 1,
$$
so that   
$$
   c_1(\sigma) 
 = \frac{1}{3}\sum_{i=1}^3 |\tau_i(\sigma)| \le \frac{7}{9}, \qquad
   c_2(\sigma)
\ge 1 - \frac{1}{3}\sum \tau_i(\wh\sigma_n) = \frac{8}{9}\,.
$$

For $n=4$,
$$
   \tau_1(\sigma) = -1, \quad
   -\frac{1}{2} \le \tau_2(\sigma) \le -(3-2\sqrt{2}), \qquad
    \frac{1}{2} \le \tau_3(\sigma) \le 4\sqrt{2}-5, \qquad
   \tau_4(\sigma) = 1,
$$
so that   
$$
   c_1(\sigma) 
\le \frac{1}{4}\sum_{i=1}^4 |\tau_i(\sigma)| >  0.78, \qquad
   \Omega_{n-1}(1,\sigma) 
\ge 1 - \frac{1}{4}\sum \tau_i(\wh\sigma_n) = 0.87\,.
$$

c$_2$) In the case $\sigma \in [\wh\sigma_n,\sigma]$, we have
$$
    \tau_i(\wh\sigma) 
 =  (1 + t)\cos\frac{(n-i)\pi}{n}  - t, \qquad 
    t \in [0,\tan^2\frac{\pi}{2n}],
$$
and, for $n=2$,
$$
   c_1(\sigma)
= \frac{1}{2}\sum_{i=1}^4 |\tau_i(\sigma)| = \frac{1+t}{2}, \qquad
   c_2(\sigma)
 = 1 - \frac{1}{2}\sum \tau_i(\sigma) = \frac{1+t}{2}\,,
$$
while for $n=3$,
$$
   c_1(\sigma)
= \frac{1}{3}\sum_{i=1}^3 |\tau_i(\sigma)| = \frac{2+t}{3}, \qquad
   c_2(\sigma)
 = 1 - \frac{1}{3}\sum \tau_i(\sigma) = \frac{2+2t}{3}\,,
$$
whereas for $n=4$,
$$
   c_1(\sigma)
 = \frac{1}{4}\sum_{i=1}^4 |\tau_i(\sigma)| 
 = \frac{\sqrt{2}+1}{4}(1+t), \qquad
   c_2(\sigma) 
 = 1 - \frac{1}{4}\sum \tau_i(\sigma) = \frac{3+3t}{4}\,.
$$


\section{Lower bounds for $B_{n,k}$}


In Lemma \ref{m>B}, we proved that
$$
   m_k(1, I_s) \ge B_{n,k}\,,
$$
where $B_{n,k}$ is the best constant in the 
Landau-Kolmogorov inequality on the half-line subject to normalization
as given below:
$$
     B_{n,k} 
 = \sup \{|f^{(k)}(1)|:\; \|f\|_{[-\infty,1]} = \|T_n\|, \;
     \|f^{(n)}\|_{[-\infty,1]} = \|T_n^{(n)}\| \}\,.
$$
So, any lower bound for $B_{n,k}$ serves as a lower bound for $m_k(1,I_s)$.

If $g$ is an arbitrary function from $W^n_\infty[-\infty,1]$, then its linear 
transfromation
$$
     f(x) 
:= \frac{\|T_n\|}{\|g\|} 
  g\left(x\cdot\Big(\frac{\|g\|}{\|T_n\|}
      \frac{\|T_n^{(n)}\|}{\|g^{(n)}\|}\Big)^{1/n}\right)
$$
is a properly normalized function, and 
$$
   B_{n,k} \ge \sup_{f} |f^{(k)}(1)|
=  \sup_g \frac{ |g^{(k)}(1)| }{ \|g\|^{1-k/n} \|g^{(n)}\|^{k/n} }
   \|T_n\|^{1-k/n} \|T_n^{(n)}\|^{k/n} 
=:  \gamma_{n,k} T_n^{(k)}(1)\,,
$$
where
$$
   \gamma_{n,k} = C_{n,k}/T_{n,k}, \qquad
   C_{n,k} 
:= \sup \frac{ |g^{(k)}(1)| }{ \|g\|^{1-k/n} \|g^{(n)}\|^{k/n} }\,, \qquad
   T_{n,k} 
:= \frac{ |T_n^{(k)}(1)| }{ \|T_n\|^{1-k/n} \|T_n^{(n)}\|^{k/n} }\,.
$$
The constant $C_{n,k}$ is the best constant 
in the LK-inequality on the half-line in the homogeneous form
$$
    |g^{(k)}(1)| \le C_{n,k} \|g\|^{1-k/n} \|g^{(n)}\|^{k/n}\,.
$$
Stechkin proved that
$$
      C_{n,k} 
\ge \frac{k!}{(2k)!}\Big(\frac{(2n)!}{n!}\Big)^{k/n} \quad
  \mbox{and}\quad    C_{n,k} 
\ge \frac{(2n!)^{1-k/n}}{(n-k)!}\,,
$$
whichever is preferrable. He also showed that 
$$
     a \Big(\frac{n}{p}\Big)^p
\le C_{n,k} \le T_{n,k} \le 
     A \Big(\frac{2n}{p}\Big)^p\,, \qquad p = \min (k,n-k)\,.
$$

\begin{lemma}
We have 
$$
    B_{n,k} \ge \gamma_{n,k} T_n^{(k)}(1),
$$
where
$$
   \gamma_{n,k} \ge (2/e)^{2k}
$$
\end{lemma}

\proof
We have
$$
    C_{n,k} = \frac{k!}{(2k)!} \Big(\frac{(2n)!}{n!}\Big)^{k/n}\,,\qquad
    T_{n,k} 
 = \frac{2^k k!}{(2k)!} 
   \frac{n^2(n^2-1^2)\cdots(n^2-(k-1)^2)}{(2^{n-1} n!)^{k/n} }\,,
$$
so that
\ba
      \gamma_{n,k} 
\ge C_{n,k}/T_{n,k} 
&=& \frac{2^{-k/n}}{n^2(n^2-1^2)\cdots(n^2-(k-1)^2)} (2n)!^{k/n} 
\lb{est1} \\
&>&  \frac{n^{2k}}{n^2(n^2-1^2)\cdots(n^2-(k-1)^2)}(2/e)^{2k} 
 > (2/e)^{2k}\,, \lb{est2}
\ea
where we used
$$
   (2n)!^{k/n} 
> \left(\sqrt{4\pi n} (2n/e)^{2n}\right)^{k/n}
> 2^{k/n} n^{2k} (2/e)^{2k}\,.
$$

\begin{lemma} 
For $n \le 15$, and $1 \le k \le n-1$, we have
$$
   C_{n,k} > \frac{ T_{n+m}^{(k)}(1)}{ T_{n+m}^{(n)}(1)^{k/n} }\,,
$$
where 
$$
   m = 1, \quad 3 \le n \le 6, \qquad
   m = 2, \quad 7 \le n \le 10, \qquad
   m = 3, \quad 11 \le n \le 14\,.
$$
\end{lemma}

\proof
For $x \in [-1,1]$, consider the function
$$
     g(x) := g_{n,m}(x) := \phi(x) T_{n+m}(x), \qquad 
     \phi(x) = c_n\,\int_{-1}^x (1-t^2)^{n}\,dt, \quad \phi(1) = 1\,,
$$
where the last equality defines the constant $c_n$. We extend it 
to the half-line $[-\infty,1$] by setting $g_{n,m}(x) = 0$ for $x < -1$. 
Then
$$
     g \in W^n_\infty[-\infty,1], \qquad 
     g^{(k)}(1) = T_{n+m}^{(k)}(1), \quad k = 1,\ldots,n\,,
$$
and 
$$ 
   C_{n,k} 
\ge \frac{ |g^{(k)}(1)| }{ \|g\|^{1-k/n} \|g^{(n)}\|^{k/n} }
= \frac{ |g^{(k)}(1)| }{ \|g^{(n)}\|_{[-1,1]}^{k/n} }\,.
$$
So, we are done once we prove that $\|g^{(n)}\|_{[-1,1]} = g^{(n)}(1)$. 
The latter is proved numerically: the graph of the function 
$g^{(n)} = g_{n,m}^{(n)}$ (provided by MAPLE) shows that, on $[-1,1]$,
for the values $n$ and $m$ given above, it 
attains its maximum at $x = 1$.

\begin{corollary}
We have 
$$
     m_k(1,I_s) > \gamma_{n,k} T_n^{(k)}(1)
$$
where 
\be
    \gamma_{n,k} 
= \frac{T_{n+m}^{(k)}(1)}{T_n^{(k)}(1)}
  \left(\frac{ T_n^{(n)}(1)}{T_{n+m}^{(n)}(1)}\right)^{k/n}\,.
\lb{gT}
\ee
\end{corollary}



\section{Proof of Theorem \ref{2} \lb{T2}}


1) For $k=1$, we have the inequality
$$
   m_1(1,I_s) \ge B_{n,1} = \gamma_{n,1} T_n'(1),
$$
where, by \rf[est1]-\rf[est2], 
$$
   \gamma_{n,1} = \frac{2^{-1/n}}{n^2}(2n)!^{1/n}
> (2/e)^2 > 0.541, \qquad \gamma_{3,1} > 0.79.
$$
We proved in \rf[m1] that
$$
   m_1(x,\sigma_n) \le \frac{1}{2}\,T_n'(1)\,, \qquad x \in [0,\omega_1]\,,
$$
and we also have 
$$
    m_1^*(x,\sigma_n) \le \alpha_{n,1} T_n'(1), \qquad x \in [\omega_1,1]\,,
$$
where
$$
    \alpha_{n,1} = \frac{1}{3} \left( \frac{2(2n-2)}{n-2}\right)^{1/n} 
\le \alpha_{4,1} = \frac{1}{3}6^{1/4} < 0.522, \quad n \ge 4, 
    \qquad \alpha_{3,1} = 2/3\,.
$$

2) For $k=2$, we have
$$
   m_2(1,I_s) \ge B_{n,2} \ge \gamma_{n,2} T_n''(1),
$$
where
$$
   \gamma_{n,2} = \frac{2^{-2/n}}{n^2(n^2-1)}(2n)!^{2/n} 
   > (2/e)^4 > 0.293\,.
$$
For the upper bounds, we have
$$
    m_2^*(x,\sigma_n) \le \alpha_{n,2} T_n''(1), \qquad
    \alpha_{n,2} = 0.23\left(\frac{2(2n-3)}{(n-3)}\right)^{2/n}
$$
and 
$$
   m_2(x,\sigma_n) \le \beta_{n,2} \,T_n''(1)\,,
$$
where
$$
   \beta_{n,2} 
= \frac{1}{5}\frac{8}{11}\frac{1}{(1-\sin\frac{3\pi}{2n})^2} < 0.28, 
  \quad n \ge 16, \qquad
   \beta_{n,2} 
=  \frac{3}{5}, \quad n < 16\,.
$$
We put all the values in the table.
$$
\begin{array}{*{4}{|c}|}
\hline
             & n =4  & 5 \le n \le 15 & n \ge 16  \\ \hline
\alpha_{n,2} & 0.72  &  \le 0.50      & \le 0.277  \\ \hline
\beta_{n,2}  &  -    &  \le 0.60      & \le 0.288  \\ \hline
\gamma_{n,2} & 0.79  &  \ge 0.63      & \ge 0.293  \\ \hline
\end{array}
$$

\newpage

\section{Proof of Theorem \ref{3}}


The values of $\gamma_{n,k}$ in \rf[gT]:
$$
\begin{array}{|c|*{12}{c|}}
\hline
 k/n &  4     & 5      & 6      & 7      & 8      & 9      & 10  
     & 11     & 12     & 13      & 14     & 15 \\ \hline
   1 & 0.87   & 0.87   & 0.87   & 0.82   & 0.82   & 0.82 & 0.82  
     & 0.79   & 0.79   & 0.79   & 0.80   & 0.80 \\ \hline 
   2 & \c 0.79   & 0.77   & 0.77   & 0.67   & 0.67   & 0.68 & 0.68  
     & 0.63   & 0.63   & 0.64   & 0.64   & 0.65 \\ \hline 
   3 &        & \c 0.72   & 0.70   & 0.57   & 0.57   & 0.57 & 0.57  
     & 0.50   & 0.51   & 0.51   & 0.52   & 0.52 \\ \hline
   4 &        &        & \c 0.66   & \c 0.51   & 0.49   & 0.49 & 0.49  
     & 0.41   & 0.41   & 0.42   & 0.42   & 0.43 \\ \hline 
   5 &        &        &        & 0.50   & \c 0.45   & 0.43 & 0.43  
     & 0.34   & 0.34   & 0.34   & 0.35   & 0.35 \\ \hline 
   6 &        &        &        &        & 0.46   & \c 0.41 & \c 0.39  
     & \c 0.30   & 0.29   & 0.29   & 0.29   & 0.29 \\ \hline 
   7 &        &        &        &        &        & 0.43 & 0.38
     & 0.27   & \c 0.26   & \c 0.25   & 0.25   & 0.25 \\ \hline 
   8 &        &        &        &        &        &      & 0.41
     & 0.27   & 0.25   & 0.23   & \c 0.22   & \c 0.22 \\ \hline 
   9 &        &        &        &        &        &        & 
     & 0.31   & 0.25   & 0.22   & 0.21   & 0.20 \\ \hline 
  10 & & & & & & & 
     &        & 0.29   & 0.23   & 0.20   & 0.19 \\ \hline 
  11 & & & & & & &  
     &        &        & 0.28   & 0.22   & 0.19 \\ \hline 
  12 & & & & & & & 
     &        &        &        & 0.27   & 0.20 \\ \hline 
  13 & & & & & & & 
     &        &        &        &        & 0.26 \\ \hline 
\end{array}
$$

\medskip
The values of 
$$
    \alpha_{n,k} 
= \frac{1}{k+1}\frac{n-1}{n-1+k}
  \left(\frac{2(2n-(k+1))}{n-(k+1)}\right)^{k/n}
$$
$$
\begin{array}{*{14}{|c}|}
\hline
 k/n &  4     & 5      & 6      & 7      & 8      & 9      & 10  
     & 11     & 12     & 13     & 14     & 15 \\ \hline
   1 & 0.58   & 0.55   & 0.54   & 0.53   & 0.53   & 0.52   & 0.52
     & 0.52   & 0.51   & 0.51   & 0.51   & 0.51 \\ \hline  
   2 & \c 0.63   & 0.48   & 0.43   & 0.40   & 0.39   & 0.38   & 0.37
     & 0.36   & 0.36   & 0.36   & 0.35   & 0.35 \\ \hline 
   3 &        & \c 0.63   & 0.44   & 0.37   & 0.34   & 0.32   & 0.30
     & 0.30   & 0.29   & 0.28   & 0.28   & 0.28 \\ \hline 
   4 &        &        & \c 0.64   & \c 0.42   & 0.34   & 0.30   & 0.28
     & 0.26   & 0.25   & 0.24   & 0.24   & 0.23 \\ \hline 
   5 &        &        &        & 0.65   & \c 0.40   & 0.32   & 0.28
     & 0.25   & 0.24   & 0.22   & 0.22   & 0.21 \\ \hline 
   6 &        &        &        &        & 0.67   & \c 0.40   & \c 0.31
     & \c 0.26   & 0.24   & 0.22   & 0.21   & 0.20 \\ \hline 
   7 &        &        &        &        &        & 0.68   & 0.40
     & 0.30   & \c 0.25   & \c 0.22   & 0.20   & 0.19 \\ \hline 
   8 &        &        &        &        &        &        & 0.69
     & 0.39   & 0.29   & 0.24   & \c 0.21   & \c 0.19 \\ \hline 
   9 &        &        &        &        &        &        &
     & 0.70   & 0.39   & 0.29   & 0.24   & 0.21 \\ \hline 
  10 & & & & & & & 
     &        & 0.71   & 0.39   & 0.29   & 0.23 \\ \hline 
  11 & & & & & & &  
     &        &        & 0.72   & 0.39   & 0.28 \\ \hline 
  12 & & & & & & & 
     &        &        &        & 0.73   & 0.39 \\ \hline 
  13 & & & & & & & 
     &        &        &        &        & 0.74 \\ \hline 
\end{array}
$$

\medskip
It is readily seen that, for the values of $k\backslash n$
above and on the shadowed cells, we have 
$$
    \alpha_{n,k} \le \gamma_{n,k}\,.
$$





%


\end{document}